\documentclass{amsart}
\usepackage{amssymb}
\usepackage{graphicx}

\def\alp{\alpha}

\def\dd{\displaystyle}

\def\gam{\gamma}
\def\Gam{\Gamma}

\def\mtline#1{\hbox to#1{\hrulefill}}

\def\noi{\noindent}

\def\ome{\omega}
\def\Ome{\Omega}

\def\sig{\sigma}

\def\what{\widehat}
\def\wtit{\widetilde}

\def\cA{{\mathcal A}}

\def\cG{{\mathcal G}}

\def\cK{{\mathcal K}}

\def\cM{{\mathcal M}}

\def\cO{{\mathcal O}}

\def\cR{{\mathcal R}}
\def\cS{{\mathcal S}}

\def\cU{{\mathcal U}}

\def\cW{{\mathcal W}}

\def\maA{{\mathcal A}}

\def\maG{{\mathcal G}}

\def\maR{{\mathcal R}}
\def\maS{{\mathcal S}}

\def\maU{{\mathcal U}}

\def\D{{\mathbb D}}

\def\N{{\mathbb N}}

\def\R{{\mathbb R}}
\def\S{{\mathbb S}}

\def\Z{{\mathbb Z}}

\def\wtL{\widetilde{L}}

\def\fd{\mathfrak d}

\def\fg{\mathfrak g}

\def\fl{\mathfrak l}

\usepackage{secteqn}
\usepackage{mathabx}

\usepackage{color}
\definecolor{purple}{cmyk}{.33,1,0,.4}
\definecolor{m}{rgb}{1,0.1,1}
\definecolor{green}{cmyk}{1,0,1,0}
\definecolor{test}{rgb}{1,1,1}
\definecolor{cmyk}{cmyk}{0,1,1,0}

\newtheorem{Equation}{}[section]

\newtheorem{definition}[Equation]{Definition}
\newtheorem{example}[Equation]{Example}

\newtheorem{lemma}[Equation]{Lemma}
\newtheorem{proposition}[Equation]{Proposition}

\newtheorem{remark}[Equation]{Remark}
\newtheorem{remarks}[Equation]{Remarks}
\newtheorem{theorem}[Equation]{Theorem}

\newtheorem{coro}[Equation]{Corollary}
\newtheorem{defi}[Equation]{Definition}

\newtheorem{prop}[Equation]{Proposition}

\def\AS{\operatorname{AS}}
\def\ch{\operatorname{ch}}

\def\diag{\operatorname{diag}}

\def\dim{\operatorname{dim}}

\def\Dom{\operatorname{Dom}}

\def\End{\operatorname{End}}

\def\cohom{\operatorname{H}}

\def\Id{\operatorname{I}}

\def\Ind{\operatorname{Ind}}

\def\oK{\operatorname{K}}

\def\Ran{\operatorname{Ran}}

\def\Tr{\operatorname{Tr}}

\def\tr{\operatorname{tr}}

\def\hTr{\operatorname{\mathfrak{Tr}}}

\def\vol{\operatorname{vol}}

\def\cA{{\mathcal A}}

\def\maA{{\mathcal A}}

\def\alp{\alpha}

\def\gam{\gamma}
\def\Gam{\Gamma}

\def\ep{\epsilon}

\def\lam{\lambda}
\def\Lam{\Lambda}

\def\sig{\sigma}

\def\ome{\omega}
\def\Ome{\Omega}

\def\lan{\langle}
\def\ran{\rangle}

\def\dd{\displaystyle}

\def\ovln{\overline}
\def\pa{\partial}

\def\ms{\medskip}

\def\ssm{\smallsetminus}

\def\what{\widehat}
\def\wtit{\widetilde}

\begin{document}



\title[Higher relative  theorems]
{ Higher relative index theorems for foliations \\ \today}


\author{Moulay Tahar Benameur}
\address{Institut Montpellierain Alexander Grothendieck, UMR 5149 du CNRS, Universit\'e de Montpellier}
\email{moulay.benameur@umontpellier.fr}

\author[James L. Heitsch]{James L.  Heitsch}
\address{Mathematics, Statistics, and Computer Science, University of Illinois at Chicago} 
\email{heitsch@uic.edu}

\thanks{\hspace{-0.5cm}
}
\begin{abstract}  In this paper we solve the general case of the cohomological relative index problem for foliations of non-compact manifolds. In particular, we significantly generalize the groundbreaking results of Gromov and Lawson, \cite{GL83}, to Dirac operators defined along the leaves of  foliations of non-compact complete Riemannian manifolds, by involving all the terms of the Connes-Chern character, especially the higher order terms in  Haefliger cohomology. The zero-th order term corresponding to holonomy invariant measures was carried out in  \cite{BH21} and becomes a special case of our main results here.  
In particular, for two leafwise Dirac operators on two foliated manifolds which agree near infinity, we define a relative topological index and the Connes-Chern character of a relative analytic index, both being in relative Haefliger cohomology.  We show that these are equal.  This invariant can be paired with closed holonomy invariant currents (which agree near infinity) to produce higher relative scalar invariants. 
When we relate these invariants to the leafwise index bundles, we restrict to Riemannian foliations on manifolds of sub-exponential growth. This allows us {to prove a higher relative index bundle theorem}, extending the classical index bundle theorem of \cite{BH08}.  Finally, we construct examples of foliations and use these invariants to prove that their spaces of leafwise positive scalar curvature metrics have infinitely many path-connected components, completely new results which  are not available from \cite{BH21}. In particular, these results confirm  the well-known idea that important geometric information of foliations is embodied in the higher terms of the  $\what{A}$ genus. 

\end{abstract} 

\maketitle

\tableofcontents

\section{Introduction}

In this paper we continue our program of extending the groundbreaking relative index {theorems of Gromov-Lawson, especially Theorem 4.18, \cite{GL83},  to Dirac operators} defined along the leaves of foliations of non-compact complete Riemannian manifolds.  Their result{s have} played a fundamental role in the development and understanding of the existence and non-existence of metrics with positive scalar curvature (PSC), as well as the structure of  spaces of such metrics.  It is {an essential tool} for the extension of results for compact manifolds to non-compact manifolds.  

\ms

In \cite{BH21}, we extended the Gromov-Lawson theorem to foliations admitting invariant transverse measures, and  crucial requirements for the applications were that the foliation admits a holonomy invariant measure, and that the measured $\what{A}$ genus of the foliation be non-zero.  In this paper, we dispense with both these requirements and completely solve the general case.  We obtain results for all  the terms of the Atiyah-Singer characteristic forms associated with the Dirac operators, especially the higher order terms of the  Connes-Chern character of the relative analytic index, as well as the higher order terms of the Connes-Chern characters of their ``index bundles".  
We also construct a large collection of spin foliations, with trivial zero-th order Haefliger $\what{A}$ genus, whose spaces of leafwise PSC metrics have infinitely many path connected components.  In particular, these results confirm the idea that the higher order terms of the $\what{A}$ genus  carry important geometric information. 
 
\ms

As in \cite{BH21}, our work is in the spirit of the transition from the Atiyah-Singer index theorem, \cite{ASIII}, to Connes'  index theorem for foliations, \cite{C79, C81, CS84}.  In order to overcome the problems of dealing with non-compact manifolds, we assume that our objects have bounded geometry.  Our higher relative index theorem then provides the expected formula in an appropriate  relative Haefliger cohomology  for pairs of foliations which are isomorphic near infinity, equating the higher relative analytical index constructed out of parametrices with the higher relative $A$-hat forms. When the foliations are top dimensional, we recover the Gromov-Lawson theory \cite{GL83, LM89}. When the foliated manifolds are compact (without boundary), we recover the cohomological version of the Connes-Skandalis index theorem \cite{CS84}, as  developed in \cite{BH04} using Haefliger cohomology. When the foliations are not top-dimensional, any pair of Haefliger transverse currents which are compatible near infinity lead to scalar higher relative index formulae. We thus recover the results of our previous paper \cite{BH21} by pairing our higher relative index formula with a compatible pair of holonomy invariant transverse measures.

\ms

 As is well known and already observed for closed foliated manifolds, see for instance \cite{BH08, BHW14, HL99}, despite the top-dimensional case, further conditions are required to relate the higher analytic index of leafwise Dirac operators to their spectral index, say the Connes-Chern characters of the leafwise projections to their kernels, the so-called index bundle. The examples in \cite{BHW14} show that such restrictions are necessary. 
Assuming, as in \cite{BH08}, that the spectral projections of the leafwise Dirac operators are sufficiently sparse near zero and that the foliations are Riemannian, we prove our next higher index theorem which now involves the relative spectral index. This theorem holds only in the absolute Haefliger cohomologies since the pair of index bundles is in general  not compatible near infinity.  This incompatibility can prevent the pairing of the index bundles with compatible near infinity Haefliger currents from being well defined. Finally, we show that when the ambient manifolds have sub-exponential growth, such pairings are miraculously well defined as soon as the Dirac operators are invertible near infinity, and they equal the pairing with the higher relative $A$-hat forms.  The invertibility near infinity is the usual Gromov-Lawson condition involving the zero-th order term of the Bochner formula.  It occurs for instance when the foliations are spin with leafwise PSC near infinity, compare with \cite{GL83}.  
\ms

Notational details are given in the next section.

\ms

Denote by $(M, F)$ a foliated manifold where $M$ is a non-compact complete Riemannian manifold and $F$ is an oriented foliation (with the induced metric) of $M$.  We assume that both $M$ and $F$ are of bounded geometry and that the holonomy groupoid of $F$ is Hausdorff.  We will sometimes assume that $F$ is Riemannian, and when we do, we will explicitly point it out in the text.  The general case will be addressed in \cite{BH23}. We assume that we have a Clifford bundle $E_M \to M$ over the Clifford algebra of the co-tangent bundle to $F$, along with a Hermitian connection $\nabla^{F,E}$ compatible with Clifford multiplication.  This determines a leafwise generalized Dirac operator, denoted  $D_F$.  We assume that we have a second foliated manifold $(M',F')$ with the same structures.  We further assume that there are compact subspaces $\cK_M = M\smallsetminus V_M$ and $\cK'_{M'} = M'\smallsetminus V'_{M'}$ so that the situations on $V_M$ and $V'_{M'}$ are identical via a smooth isometry $\varphi$. These are the usual Gromov-Lawson relative data.  {{Note that in our case, the ``bad set" restricted to a leaf need {\em not} be compact as in the Gromov-Lawson case.  Only the global aggregate of all such leafwise subsets needs to be compact as a subset of $M$.}}

\ms

In \cite{BH21}, we worked on the ambient manifolds $M$ and $M'$.  Here we work on their holonomy groupoids $\cG$ and $\cG'$, with their canonical foliations $F_s$ and $F'_s$, as we did in \cite{BH08}.  We lift everything to $\cG$ using the {range} map $r:\cG \to M$, which is a covering map from the leaves of $F_s$ to those of $F$, and similarly for $M'$. In particular, we have the $\cG$ invariant leafwise Dirac operator $D$ for the foliation $F_s$, and similarly $D'$ for $F'_s$.

\ms

Recall that for a good cover $\cU = \{(U_i,T_i) \}$ of $M$, \cite{HL90}, by foliation charts $U_i$ with local complete transversals $T_i\subset U_i$, the Haefliger forms associated to $F$ are the bounded smooth differential forms on $\amalg\,T_i$ which have compact support in each $T_i$, modulo forms minus their holonomy images.  The (absolute) Haefliger cohomology of $F$, denoted $H^*_c(M/F)$, is then the associated de Rham cohomology, and is independent of the choice of good cover, \cite{Ha80}.  Also recall that there is an {\em{integration over the leaves}}  map from forms on $M$ to Haefliger forms, denoted $\dd \int_F$, which induces a map on cohomology.  For the foliation given by the fibers of a bundle $M \to B$, the Haefliger cohomology reduces to the cohomology of the base and $\dd \int_F$ is the classical integration over the fibers map.  See again \cite{Ha80} for more details.

\ms

The receptacle for our relative index formulae will be a relative version of Haefliger cohomology that we denote by $H^*_c (M/F, M'/F';\varphi)$. This is the cohomology of pairs of Haefliger forms which agree near infinity (that is, on $T_i$ far enough away from $\cK_M$ and similarly for the $T_i'$), again modulo pairs of forms minus their holonomy images which also agree near infinity.
\ms

Denote by $AS(D_F)$ the Atiyah-Singer characteristic differential form, associated with the above $\varphi$-compatible data, for $D_F$, and similarly for $D_{F'}$.  These differential forms agree near infinity on $M$ and $M'$.  The relative $A$-hat genus of  {{the compatible pair $(D, D')$, alternatively called the relative topological  index,} is
$$
\Ind_t (D, D') \; = \; \left[\int_F \AS(D_F),\int_{F'}AS(D'_{F'}))\right]  \; \in \; H^*_c (M/F, M'/F';\varphi).
$$
 Using parametrices, we define a relative  analytical index class $\Ind_a (D, D')$ in the appropriate $K$-theory group, and its Connes-Chern character,
$$
{{\ch}}(\Ind_a (D, D')) \;  \in \; H^*_c (M/F, M'/F';\varphi).
$$  
Our first result is

\ms

{\bf{Theorem \ref{classicalindex} }}
{\em {For the $(M,F)$, $(M',F')$, $D$ and $D'$ as above,  
$$
\ch(\Ind_a (D, D')) = \Ind_t (D, D')\; \text{ in }H^*_c (M/F, M'/F';\varphi).
$$}}
So, pairing with any compatible near infinity pair $(C, C')$ of closed Haefliger currents yields a (higher) scalar  relative index formula.  Such pairings will be denoted $\langle \cdot, \cdot \rangle$, e.g. $\langle \ch(\Ind_a (D, D')), (C, C')\rangle$.

\ms

An important application of this theorem is to pairs of ``reflective" foliations, which we consider in Section 6.  They can be ``cut and pasted" to get a compact foliated manifold $\what{M}$, with the foliation $\what{F}$ and operator $\what{D}_{\what{F}}$.   Given $C$ and $C'$ as above, denote by $\what{C}$ the current they determine on $\what{M}$.  Then we have the following extension of {the Gromov-Lawson Relative Index Theorem, see \cite{GL83}}, 
which is most useful in Section \ref{apps}, where we construct our examples.  

\ms

\noi {\bf{Theorem \ref{ExRelInd2}} } {\em
Suppose that  $F$ (and so also $F'$)  is reflective. Then
$$
\langle \ch ( \Ind_a(D, D'), (C,C') \rangle  \: =  \: \langle \ch ( \Ind_a(\what{D}_{\what{F}})), \what{C} \rangle. 
$$
}

The RHS of this index formula can then be computed using the classical higher cohomological index theorem for foliations of closed manifolds \cite{C94, BH04}. For top dimensional foliations, say when $TF=TM$ and $TF'=TM'$, the previous two theorems reduce to the classical Gromov-Lawson relative index theorems.

\ms

Despite the top dimensional case, it is well known that  the higher index is not easily related with the so-called index bundle, i.e. the Chern character of the ``kernel minus cokernel superbundle''.   Constraints on the spectral distributions, as well as on the geometry near infinity are necessary, see for instance \cite{BHW14}. Denote by $P_0$ the leafwise spectral projection to the kernel of $D^{2}$.  In general $P_0$ is not transversely smooth (although it is always leafwise smooth), and if not, we cannot even define its Connes-Chern character in our Haefliger cohomology without perturbing the operator.  There are though interesting classes of foliations and leafwise Dirac-type operators whose kernel superbundle $P_0$ is transversely smooth, and in this case, we get a {well defined spectral index class
$$
\ch(P_0)  \;  \in \; H_c^*(M/F),
$$
and similarly for $P_0'$, see \cite{BH08}.

\ms

Denote by $P_{(0,\ep)}$ the leafwise spectral projection for $D^{2}$ for the interval $(0,\ep)$. 
The Novikov-Shubin invariants $NS(D)$ of $D$ are a measure of the density of the image of $P_{(0,\ep)}$.  The larger $NS(D)$ is, the sparser the image of $P_{(0,\ep)}$ is as $\ep \to 0$.

\ms

We also have the natural map $(\pi\times \pi'): H^*_c (M/F, M'/F';\varphi) \to H_c^*(M/F)\times H_c^*(M'/F')$,
and with it the Riemannian Foliation Relative Index Bundle Theorem. 

\ms

\noi
{\bf{Theorem} \ref{Secondhalf1}} {\em 
Fix $0 \leq \ell \leq q/2$, where $q$ is  the codimension of $F$ and $F'$.  Assume that:
\begin{itemize} 
\item the foliations $F$ and $F'$ are Riemannian;
\item the leafwise operators $P_0$, $P_0'$, $P_{(0,\ep)}$  and $P_{(0,\ep)}'$ (for $\epsilon$ sufficiently small) 
are transversely smooth;
\item  $NS(D)$ and $NS(D')$ are greater than {{$\ell$}}.
\end{itemize}
Then,  
for $0 \leq k \leq \ell$,  we have in $H_c^{2k}(M/F)\times H_c^{2k}(M'/F')$
$$
(\pi\times \pi') \ch^k(\Ind_a(D, D'))\, = \,   (\ch^k(\Ind_a(D)), \ch^k(\Ind_a(D'))) \, = \, 
(\ch^k(P_{0}), \ch^k(P'_{0})).
$$}

For Riemannian foliations, important examples of compatible near infinity pairs of closed Haefliger currents are given by closed bounded holonomy invariant transverse differential forms $\ome$ on $M$ and $\ome'$ on  $M'$ which agree near infinity.  These determine closed bounded Haefliger forms on $T$, denoted $\ome_T$ and $\ome'_{T'}$ which agree near infinity.  Denote by $dx$ the global volume form on $M$.   

\ms

We then have the Higher Relative Index Pairing Theorem. 

\smallskip

\noi
{\bf{Theorem \ref{Secondhalf2} }} {\em 
In addition to the assumptions in Theorem \ref{Secondhalf1}, assume that for $\epsilon$ sufficiently small,  $\dd \int_M \tr (P_{[0,\ep)}) dx  <  \infty$ and $\dd \int_{M'} \tr (P'_{[0,\ep)}) dx  <  \infty$, and that
$M$, and so also $M'$, has sub-exponential growth. Then, for any $\ome  \in C^{\infty}(\wedge^{q-2k} \nu^*)$ and $\ome' \in C^{\infty}(\wedge^{q-2k}  {\nu'}^*)$ ($0 \leq k \leq \ell$) as above, 
$$
\dd \int_T {{\ch}}(P_0) \wedge \ome_T  \text{ and }\dd \int_{T'} {{\ch}}(P'_0) \wedge \ome'_{T'}  \text{  are well defined complex numbers,}
$$
and
$$ 
\int_T {{\ch}}(P_0) \wedge \ome_T \, - \, \int_{T'} {{\ch}}(P'_0) \wedge \ome'_{T'}  \, = \
\langle \left[\int_F \AS(D_F), (\int_{F'}AS(D'_{F'})\right],   [\ome_T,\ome'_{T'}]\rangle. 
$$ }

\ms

In Section \ref{TIMAFTS},  we  show that the finite integral assumptions in Theorem \ref{Secondhalf2} are satisfied when $D_F$ (and hence also $D_{F'}$) is invertible near infinity, i.e. when the zeroth order differential operator $\maR^E_F$ in the associated Bochner Identity 
 $$
 D_F^2=\nabla^*\nabla + \maR^E_F, 
 $$ 
 is uniformly positive near infinity on $M$.  The sub-exponential growth condition can be extended to exponential growth provided it is not too robust.  See Remark \ref{growth}. 

\ms

For a single foliated manifold with a pair of compatible near infinity leafwise Dirac operators, we have the following generalization of a classical result of Gromov-Lawson \cite{GL83}, compare with Theorem 6.5 in \cite{LM89}.

\ms

\noi 
{\bf{Theorem \ref{IMCOR2}} }  {\em
Suppose that $E$ and $E'$ are two Clifford bundles over the foliated manifold $(M,F)$, which are isomorphic off the compact subset $\cK_M$, with associated twisted Dirac operators $D$ and $D'$.  Let $\omega$ be a bounded closed holonomy invariant transverse form (or Haefliger current) of degree $\ell\leq q$.  Suppose that
\begin{itemize}
\item $M$ has sub-exponential growth, and $F$ is Riemannian;
\item the leafwise operators $P_0$, $P_0'$, $P_{(0,\ep)}$  and $P_{(0,\ep)}'$ (for $\epsilon$ sufficiently small)  are transversely smooth;
\item $\min (NS(D), NS(D'))$ is greater than $\ell$;
\item $\maR^E_F$, and hence also $\maR^{E'}_F$,  is uniformly positive near infinity.   
\end{itemize}
 Then, since $\ch(E)  = \ch(E')$ off $\cK_M$,
$$
\int_{\cK_M} (\AS (D_F)(\ch(E) - \ch(E')) \wedge \ome 
\, = \,  
 \int_T \left(\ch(P_0) \, - \, {{\ch}}(P'_0)\right) \wedge \ome_T.
$$  }

In the reflective case, again more constraints are necessary to obtain the link with the index bundle, and we have the following.

\smallskip

\noi{\bf{Theorem \ref{ExRelInd3} }} {\em 
Suppose that  $F$ (so also $F'$) is reflective. Suppose furthermore that $\what F$ is Riemannian and that $\what{P}_0$ and $\what{P}_{(0, \ep)}$ are transversely smooth and the Novikov-Shubin invariants of $\what{D}_{\what{F}}$ are greater than $\ell$, for some $0\leq \ell \leq q/2$.  Then for any $2\ell$ homogeneous $\varphi$-compatible $(\ome, \ome')$ as above,
$$ 
\langle \ch( \Ind_a(D,D')), {{[\ome_T,\ome'_{T'}]}}\rangle  =  \langle(\ch(\what{P}_{0}),\what{\ome}_{\what{T}}\rangle.
$$
Moreover, if we impose the assumptions of Theorem \ref{Secondhalf2}, then 
$$
\langle (\ch(P_0), \ch(P_0)),  (\ome_T, \ome'_{T'})  \rangle \; = \; \langle(\ch (\what{P}_{0}),\what{\ome}_{\what{T}}\rangle.
$$   }

In Section \ref{apps}, we consider foliations which admit positive scalar curvature (PSC) leafwise metrics.  Given such a foliation, we associate to any pair $(g_0, g_1)$ of such metrics, an invariant living in Haefliger cohomology, which provides an obstruction for the leafwise path connected equivalence of $g_0$ and $g_1$. This precisely generalizes the classical Gromov-Lawson invariant.  Finally, we construct a large collection of spin foliations whose space of leafwise PSC metrics has infinitely many path connected components.

\ms

\noi
{\em Acknowledgements.} 
MTB thanks the french National Research Agency for support via the project ANR-14-CE25-0012-01 (SINGSTAR).
JLH thanks the Simons Foundation for a Mathematics and Physical Sciences-Collaboration Grant for Mathematicians, Award Number 632868.

\section{The Setup}\label{setup}

Denote by  $M$ a  smooth non-compact complete Riemannian manifold of dimension $n$, and by $F$  an oriented foliation (with the induced metric) of $M$ of dimension $p$, {{(until further notice, we assume that $p$ is even), }}
and codimension $q = n - p$.   The tangent and cotangent bundles of $M$ and $F$ are denoted $TM, T^*M, TF$ and $T^*F$. The normal and dual normal bundles of $F$ are denoted $\nu$ and $\nu^*$.  A leaf of $F$ is denoted by $L$.  At times, we will assume that $F$ is Riemannian, that is the metric on $M$, when restricted to $\nu$ is bundle like, so the holonomy maps of $\nu$ and  $\nu^*$ are isometries. We will consider the general case in \cite{BH23}.   

\ms

We assume that both $M$ and $F$ are of bounded geometry, that is, the injectivity radius on $M$ and on all the leaves of $F$ is bounded below, and the curvatures and all of their covariant derivatives on $M$ and on all the leaves of $F$ are bounded (the bound may depend on the order of the derivative).   

\ms

Let  ${\cU}$ be a good cover of $M$ by foliation charts as defined in \cite{HL90}.  In particular, denote by $\D^p (r)=\{x \in \R^p,  ||x||  <r\}$, and similarly for $\D^q(r)$.   An open locally finite cover  $\{(U_i, \psi_i)\}$ of $M$  by  foliation coordinate charts $\psi_i: U_i \to  \D^p(1) \times \D^q(1) \subset \R^n$  is a good cover for $F$ provided that
\begin{enumerate}
 \item 
For each $y \in \D^q(1), P_y = \psi_i^{-1}(\D^p(1) \times  \{y\})$ is contained in a leaf of $F$. $P_y$ is called a plaque of $F$.
\smallskip
\item 
If $\overline{U}_i \cap  \overline{U}_j \neq \emptyset$,  then $U_i \cap  U_j \neq \emptyset$, and $U_i \cap  U_j$ is connected.
\smallskip
\item 
Each $\psi_i$ extends to a diffeomorphism $\psi_i: V_i \to \D^p(2) \times \D^q(2)$, so that the cover $\{(V_i, \psi_i)\}$  satisfies $(1)$ and $(2)$, with $\D^p(1)$ and $\D^q(1)$ replaced by $\D^p(2)$ and $\D^q(2)$.
\smallskip
\item 
Each plaque of $V_i$  intersects at most one plaque of $V_j$  and a plaque of $U_i$ intersects a plaque of $U_j$ if and only if the corresponding plaques of $V_i$ and $V_j$ intersect.
\smallskip
\item
There are global positive upper and lower bounds on the norms of each of the derivatives of the $\psi_i$.
\end{enumerate}
Bounded geometry foliated manifolds always admit good covers. 

\ms

For each $U_i \in {\cU}$, let $T_i\subset U_i$ be a {{local complete transversal}}  (e.g. $T_i =  \psi_i^{-1}(\{0\} \times  \D^q(1))$) and set $T=\bigcup\,T_i$.  We may assume that the closures of the $T_i$ are disjoint.  Given $(U_i,T_i)$ and $(U_j,T_j)$, suppose that $\gamma_{ij\ell}:[0,1] \to M$ is a path whose image is contained in a leaf with   $\gamma_{ij\ell}(0)\in T_i$ and $\gamma_{ij\ell}(1) \in T_j$. Then 
$\gamma_{ij\ell}$  induces a local diffeomorphism $h_{\gamma_{ij\ell}}:  T_i \to T_j$, with domain $\Dom_{\gamma_{ij\ell}}$ and range $\Ran_{\gamma_{ij\ell}}$.   The  space $\cA^k_c(T)$ consists of all smooth $k$-forms on $T $ which are $C^{\infty}$ bounded and have compact support in each $T_i$.  The  Haefliger $k$-forms for $F$, denoted  $\cA^k_c(M/F)$, consists of elements in the quotient of $\cA^k_c(T)$ by the closure of the vector subspace {{$\cW$}} generated by elements of the form $\alpha_{ij\ell}-h_{\gamma_{ij\ell}}^*\alpha_{ij\ell}$ where $\alpha_{ij\ell} \in \cA^k_c(T)$ has support contained in $\Ran_{\gamma_{ij\ell}}$.  We need to take care as to what this  means.  Members of {{$\cW$}} consist of possibly infinite sums of elements of the form $\alpha_{ij\ell}-h_{\gamma_{ij\ell}}^*\alpha_{ij\ell}$, with the following restrictions: each member of  {{$\cW$}} has a bound on the leafwise length of all the $\gamma_{ij\ell}$ for that member, and each $\gamma_{ij\ell}$ occurs at most once.   Note that these conditions plus bounded geometry imply that for each member of {{$\cW$}},  there is $n \in \N$ so that the number of elements of that member having $\Dom_{\gamma_{ij\ell}}$ contained in any $T_i$ is less than $n$,  and that each $U_i$ and each $U_j$  appears at most a bounded number of times.   The projection map is denoted 
$$
[ \cdot ]: \cA^*_c(T) \,\,  \to \,\,  \cA^*_c(M/F).
$$

\ms

Denote the exterior derivative by $d_T:\cA^k_c(T)\to \cA^{k+1}_c(T)$, which  induces $d_H:\cA^k_c(M/F)\to \cA^{k+1}_c(M/F)$.  Note that $\cA^k_c(M/F)$ and $d_H$ are independent of the choice of cover $\cU$.   The cohomology $H^*_c(M/F)$ of the complex $\{\cA^*_c(M/F),d_H\}$ is the Haefliger cohomology of $F$.

\ms

Denote by  $\cA^*_u(M)$  the space of differential forms on $M$ which are  smooth and $C^{\infty}$ bounded, and denote its exterior derivative by $d_{M}$ and its cohomology by $H_u^*(M;\R)$.  As the bundle $TF$ is oriented, there is a continuous open surjective linear map, called integration over $F$,
$$
\dd \int_F:\cA^{p+k}_u(M) \to \cA^k_c(T),
$$
which commutes with the exterior derivatives.
This map is given by choosing a partition of unity $\{\phi_i\}$ subordinate to the cover $\cU$, and setting 
$\dd \int_F \omega$ to be the class of $\dd \sum_i \int_{U_i} \phi_i \omega.$  It is a standard result,  \cite{Ha80},  that the image of this differential form $\dd \Big[ \int_F \omega \Big]   \in \cA^k_c(M/F)$ is independent of the partition of unity and of the cover $\cU$.  As $\dd \int_{F}$ commutes with $d_{M}$ and $d_{H}$, it induces the map $\dd \int_F :H_u^{p+k}(M;\R) \to H^k_c(M/F)$.

\ms

Note that $\dd \int_{U_i}$ is integration over the fibers of the projection $U_i \to T_i$, and
 that each integration $\omega \to \dd \int_{U_i} \phi_i \omega$ is essentially integration over a compact fibration, so $\dd \int_F$ satisfies the Dominated Convergence Theorem on each $U_i \in \cU$. 
  
\ms

The holonomy groupoid ${\mathcal G}$ of $F$ consists of equivalence classes of paths $\gamma:[0,1]\to M$ such that the image of $\gamma$ is contained in a leaf of $F$.  Two such paths $\gamma_{1}$ and $\gamma_{2}$  are
equivalent if $\gamma_{1}(0) = \gamma_{2}(0)$, $\gamma_{1}(1) = \gamma_{2}(1)$, and the holonomy germ along them is the same.  Two classes may be composed if the first ends where the second begins, and the composition is just the juxtaposition of the two paths.  This makes $\cG$ a groupoid.  The space $\cG^{(0)}$ of units of $\cG$ consists of the equivalence classes of the constant paths, and we identify $\cG^{(0)}$ with $M$. 

\ms

The basic open sets defining the (in general non-Hausdorff) $2p +q$ dimensional manifold  structure of $\cG$ are given as follows.   Given $U_i, U_j \in \cU$ and a leafwise path $\gamma_{ij\ell}$ starting in $U_i$ and ending in $U_j$, define the graph chart $U_i \times_{\gamma_{ij\ell}}U_j$ to be the set of equivalence classes of leafwise paths starting in $U_i$ and ending in $U_j$ which are homotopic to $\gamma_{ij\ell}$ through a homotopy of leafwise paths whose end points remain in $U_i$ and $U_j$ respectively.  It is easy to see, using the holonomy map $h_{\gamma_{ij\ell}}:T_i \to T_j$  that $U_i \times_{\gamma_{ij\ell}}U_j \simeq \D^{p}(1) \times \D^{p}(1) \times \D^{q}(1)$.   

\ms

$\cG$ has the natural the maps $r,s:\cG \to M$, with $s([\gamma]) = \gamma(0)$ and  $r([\gamma]) = \gamma(1)$.  It also has has two natural foliations, $F_s$ and $F_r$,  whose leaves are the fibers of $s$ and $r$.  We will primarily use $F_s$, whose leaves are denoted $\wtit{L}_x=s^{-1}(x)$, for $x \in M$.  Note that $r:\wtit{L}_x \to L$ is the holonomy covering map.  We will assume that $\cG$ is Hausdorff, which is automatic for Riemannian foliations. 

\ms

The smooth sections of a bundle $E$ are denoted by $C^{\infty}(E)$, and those with compact support by $C_c^{\infty}(E)$. We assume that any connection or any metric on $E$, and all their derivatives, are bounded.     See \cite{Shubin} for material about bounded geometry bundles and their properties.

\ms

For a real or complex bundle $E_M \to M$,  the external tensor product bundle $E_M \boxtimes E_M^*  \to M \times M$ can be pulled back under $(s, r)$ to a smooth bundle denoted $E\boxtimes E^*$ over $\cG$.   We denote the smooth,  bounded sections $k(\gam)$ with compact support of the restriction of this bundle to {{subset $U_i \times_{\gamma_{ij\ell}}U_j \subset \cG$}} by $C^{\infty}_c(U_i \times_{\gamma_{ij\ell}}U_j,E\boxtimes E^*)$.   We extend them to all of $\cG$ by by setting $k(\gamma) = 0$ if $\gamma \notin U_i \times_{\gamma_{ij\ell}}U_j$.   
\begin{defi}\cite{BH18}\  
The  algebra  $C^\infty_{u} (E\boxtimes E^*)$ consists of smooth sections $k$  of  $E \boxtimes E^*$, called kernels,  such that $k$ is a (possibly infinite) sum $k = \sum_{ij\ell} k_{ij\ell}$,  with each $k_{ij\ell} \in C^{\infty}_c (U_i \times_{\gamma_{ij\ell}}U_j, E\boxtimes E^* )$. 
For each $k$, we require that there is a bound on the leafwise length of its $\gamma_{ij\ell}$, and that each index $ij\ell$ occurs at most once.    We further require that for each $k$, each of its derivatives in the local coordinates given by the good cover is bounded, with the bound possibly depending on the particular derivative. 
\end{defi}
The proof of Lemma 2.3 of  \cite{BH08} shows that this is indeed an algebra. 
Each  $k \in C^\infty_{u} (E\boxtimes E^*)$ defines a $\cG$-invariant leafwise smoothing operator on $C^\infty_c (E)$ in the sense of \cite{C79}, which is transversely smooth and has finite propagation.  See  \cite{Shubin} for the definition of bounded geometry smoothing operators, as well as \cite{NistorWeinsteinXu} for the groupoid version.  To see this, use the leafwise distance function $d _x(\gam, \what{\gam})$ on $\wtit{L}_x$.   This is defined as the infimum over the leafwise length $\fl(\gam\what{\gam}^{-1})$ of all paths in the class of $\gam\what{\gam}^{-1}\in \cG$.  For any bounded geometry foliation with Hausdorff holonomy groupoid, the sets $U_i \times_{\gamma_{ij\ell}}U_j$ have the property that there is a universal constant (namely the bound $C$ on the diameters of all the placques in all the $U_i \times_{\gamma_{ij\ell}}U_j$), so that for all $\gam \in U_i \times_{\gamma_{ij\ell}}U_j$, we have 
$\fl(\gam) \leq \fl(\gam_{ij\ell}) + 2C$.  Next, suppose that $k_{ij\ell} \in C^{\infty}_c (U_i \times_{\gamma_{ij\ell}}U_j, E\boxtimes E^* )$, and $\sigma \in C^\infty_c (E)$.  Then, 
$$
k_{ij\ell}(\sig)(\gam) \; = \; \int_{\wtit{L}_{s(\gam)}}  k_{ij\ell}(\gam\what{\gam}^{-1})\sig(\what{\gam}) d\what{\gam}.
$$
Now, $k_{ij\ell}(\gam\what{\gam}^{-1}) = 0$ unless  $\gam\what{\gam}^{-1} \in U_i \times_{\gamma_{ij\ell}}U_j$, that is only if 
$\fl(\gam\what{\gam}^{-1}) =  d_{s(\gam)}(\gam, \what{\gam}) \leq \fl(\gam_{ij\ell}) + 2C$, the very definition of finite propagation. The restrictions imposed on each $k_{ij\ell}$ imply that each $U_i$ and each $U_j$  appears at most a  bounded number of times, so the sum converges locally uniformly, in particular pointwise.  These restrictions on  $k$ insure that it also has bounded propagation. 

\ms

Denote by $D_F$ a generalized  leafwise Dirac operator for the even dimensional foliation $F$.   It is defined as follows.  Let  $E_M$ be a complex vector bundle over $M$ with Hermitian metric and connection, which is of bounded geometry.   Assume that the tangent bundle $TF$ is spin with a fixed spin structure.   Because  $F$ is even dimensional, the bundle of spinors along its leaves, denoted $\cS_F$ splits as $\cS_F =\cS_F^+ \oplus \cS_F^-$.  Denote by $\nabla^F$ the Levi-Civita connection on each leaf $L$ of $F$.  $\nabla^F$ induces a connection $\nabla^F$ on $\cS_F | L$, and we denote by $\nabla^{F,E}$ the tensor product connection on ${\cS_F} \otimes E_M | L$.  These data determine a smooth family $D_F = \{D_L\}$ of leafwise Dirac operators, where $D_L$ acts on sections of $\cS_F \otimes E_M | L$ as follows.  Let $X_1,\ldots, X_p$ be a local oriented orthonormal basis of $TL$, and set
$$
D_L \; = \; \sum^p_{i=1}\rho(X_i)\nabla^{F,E}_{X_i}
$$
where $\rho(X_i)$ is the Clifford action of $X_i$ on the bundle ${\cS_F} \otimes E_M | L$.  Then $D_L$ does not depend on the choice of the $X_i$, and it is an odd operator for the $\Z_{2}$ grading of $\cS_F \otimes E_M = (\cS_F^+ \otimes E_M)\oplus (\cS_F^-\otimes E_M)$.  Thus $D_F: C^{\infty}_c(\cS_F^\pm \otimes E_M) \to C^{\infty}_c(\cS_F^\mp \otimes E_M)$, and $D_F^2: C^{\infty}_c(\cS_F^\pm \otimes E_M) \to C^{\infty}_c(\cS_F^\pm \otimes E_M)$. For more on {{the}} generalized Dirac operators {{that we are using here}}, see \cite{LM89}.  

\ms

Given a leafwise operator $A$ on $\cS \otimes E\otimes \wedge \nu_s^*$, denote its leafwise Schwartz kernel by $k_A$.  Then, depending on the context {{and under appropriate assumptions on $k_A$}}, the Haefliger traces, $\Tr(A)$ and $\hTr(A)$, of $A$ are defined to be, 
$$
\Tr(A) \; = \; \int_F \tr(k_A(\overline{x},\overline{x}))dx_F \in \cA^*_c(M/F)  \quad \text{ and } \quad
\hTr(A) \; = \; \left[ \int_F \tr(k_A(\overline{x},\overline{x}))dx_F \right]\in H^*_c(M/F),
$$
where $dx_F$ is the leafwise volume form associated with the fixed orientation of the foliation $F$.
The element  $\overline{x} \in \wtL_x$ is the class of the constant path at $x \in L \subset M$. See again for instance \cite{BH04} for more details on these constructions.

\ms

Now suppose that we have the situation in Section 4 of {{the companion paper}} \cite{BH21}.  That is, we have: 
\begin{itemize}
\item foliated manifolds $(M, F)$ and $(M', F')$;   
\item  Clifford bundles $E_M \to M$ and $E_{M'} \to M'$, with Clifford compatible Hermitian connections;  
\item  leafwise Dirac operators $D_F$ and $D_{F'}$;
\item compact subspaces $\cK_M = M\smallsetminus V_M$ and $\cK'_{M'}= M'\smallsetminus V'_{M'}$;
\item an isometry $\varphi:V_M \to V'_{M'}$  with $\varphi^{-1}(F') = F$;  
\item an isomorphism $\phi :E_M\vert_{V_M} \to E'_{M'}\vert_{V'_{M'}}$, {{covering $\varphi$}}, with $\phi^*(\nabla^{F',E'} \, |_{V'_{M'}}) = \nabla^{F,E} \, |_{V_M}$. 
\end{itemize}
The pair ${{\Phi = (\phi, \varphi)}}$ is {{thus}}  a bundle morphism from $E_M | V_M$ to $E'_{M'} |V'_{M'}$. 
The well defined (since they are differential operators) restrictions of $D_F$ and $D_{F'}$ to the sections over $V_M$ and $V'_{M'}$ agree through $\Phi$, i.e.
$$
(\Phi^{-1})^*\circ D_F \circ \Phi^* \,|_{V'_{M'}} = D_{F'} \,|_{V'_{M'}}. 
$$
Such operators are called $\Phi$ compatible.  Without loss of generality, we may assume that $\cK_M$ and $\cK'_{M'}$ are the closures of open subsets of $M$ and $M'$ respectively.

\ms

Recall the following material from \cite{BH21}.  Denote by $g: M\to [0, \infty)$ {{and $g':M'\to [0, \infty)$ compatible smooth approximations to the distance functions $\fd_M(\cK_M,x)$ and $\fd_{M'}(\cK'_{M'},x')$, where $\fd_M$ and $\fd_{M'}$ are the distance functions on $M$ and $M'$. So we assume that $g$ and $g'$ are $0$ on $\cK_M$ and $\cK'_{M'}$ respectively and they satisfy $g'\circ \varphi = g$}}. Hence, for $s \geq 0$,  the open submanifolds $M(s)=\{g>s\}$ and $M'(s)= \{g'>s\}$ agree through $\varphi$, that  is $\varphi(M(s)) = M'(s)$ and $g |_{M(s)} = g'\circ \varphi  |_{M(s)}$.     For $s\geq 0$
denote by $T_s$ the set 
$$
T_s \; = \; \{T_i\subset T \; |  \; T_i  \cap M(s) \neq   \emptyset   \}, 
$$  
and similarly for $T'_s$. 
 
\ms

Suppose that 
$(\zeta,\zeta') \in {{\cW}} \times {{\cW'}}  \subset \maA_c^*(T) \times  \maA_c^*(T')$, with  $\zeta = \sum_{(\alp, \gam)} \alpha - h^*_{\gam} \alpha$ and $\zeta' = \sum_{(\alp', \gam')} \alpha' - h^*_{\gam'} \alpha'$.   For simplicity, we have dropped the subscripts.  The  vector subspace ${{\cW}} \times_{\varphi} {{\cW'}} \subset {{\cW}} \times {{\cW'}} $  consists of elements $(\zeta,\zeta')$  which are ${\varphi}$ compatible.   This means that all but a finite number of the $(\alp, \gam)$ and $(\alp', \gam')$ are paired, that is 
$$
\alpha \; = \; \varphi^*(\alpha') \quad \text{and} \quad \gamma' \; = \; \varphi\circ \gamma, \quad \text{so} \quad 
\alpha - h^*_{\gam} \alpha \; = \;  \varphi^*( \alpha' - h^*_{\gam'} \alpha').
$$
\begin{definition}
Given  $\beta  \in \maA_c^*(T)$ and $\beta'\in\maA_c^*(T')$,  the pair $(\beta,\beta') $ is ${\varphi}$-compatible if  there exists $s \geq 0$ so that  $\beta = \varphi^*(\beta') $ on $T_s$.  Set 
$$
\maA_c^* (M/F, M'/F'; \varphi) \,\, = \,\, \{ (\beta, \beta') \in \maA_c^*(T)\times \maA_c^*(T') \, | \, (\beta,\beta')  \text{ is ${\varphi}$ compatible} \} /  (\overline{{{\cW}} \times_{\varphi} {{\cW'}} }). 
$$
\end{definition} 
The de Rham differentials on $\maA_c^*(T)$ and $\maA_c^*(T')$ yield a well defined relative Haefliger complex, whose homology spaces are denoted 
$$
H^*_c (M/F, M'/F';\varphi)=\oplus_{0\leq k\leq q} H^k_c (M/F, M'/F';\varphi),
$$ 
and there are well defined graded maps,
$$
\pi: H^*_c (M/F, M'/F';\varphi) \rightarrow H^*_c (M/F) \; \text{ and } \; \pi': H^*_c (M/F, M'/F';\varphi) \rightarrow H^*_c (M'/F').
$$
{{which are induced by the projections
$$
\maA_c^* (M/F, M'/F'; \varphi)\rightarrow \maA_c^* (M/F)\; \text{ and }\; \maA_c^* (M/F, M'/F'; \varphi)\rightarrow \maA_c^* (M'/F').
$$}}
\begin{definition}\label{pairingdefn}  Suppose $(\xi, \xi') \in \maA_c^* (M/F, M'/F'; \varphi)$, and let $C$ and $C'$ be closed (bounded) $\varphi$ compatible holonomy invariant Haefliger currents.  Set 
$$  
\lan (\xi, \xi'), (C,C') \ran \,\, = \,\,  
\lim_{s \to \infty} \left( C(\xi |_{T \ssm T_s}) -    C'(\xi' |_{T' \ssm T_s'})\right).
$$
\end{definition}
\noindent
This is well defined because any representative in $(\xi, \xi')$ is ${{\varphi}}$ compatible, so the right hand side is eventually constant.  In addition, every $(\zeta, \zeta') \in {{\cW}} \times_{\varphi} {{\cW'}} $ is ${{\varphi}}$ compatible, so satisfies
$$
\lim_{s \to \infty} \left( C(\zeta |_{T \ssm T_s}) -    C'(\zeta' |_{T' \ssm T_s'})\right) \,\, = \,\,   0.
$$
To see this, recall that there is a global bound on the leafwise length of the $\gam$ and $\gam'$ in $\zeta$ and $\zeta'$.  This, and the fact that there are only finitely many unpaired  $(\alpha,\gam)$ and  $(\alpha',\gam')$, insures that for large $s$, every unpaired $(\alpha,\gam)$ will have both $\Dom_{\gamma}$ and $\Ran_{\gamma} \subset T \ssm T_s$, so $\dd C(\alpha - h^*_{\gam} \alpha) $  will be zero, and similarly for  every unpaired $(\alpha',\gam')$.  Those $(\alpha,\gam)$ and $(\alpha',\gam')$ which are paired and appear in the integration, will have  $\Dom_{\gamma}$ and/or $\Ran_{\gamma} \subset T \ssm T_s$  with corresponding $\Dom_{\gamma'}$ and/or $\Ran_{\gamma'} \subset T' \ssm T'_s$.   In both cases, their integrals will cancel.

\begin{remark}
Examples of such currents include the following.
\begin{enumerate}
\item  Invariant transverse measures $\Lam$ and $\Lam'$ on $T$ and $T'$ which are ${{\varphi}}$ compatible as in \cite{BH21}.
\item  Suppose $\ome  \in C^{\infty}(\wedge^* \nu^*)$ and $\ome' \in C^{\infty}(\wedge^* {\nu'}^*)$ are closed holonomy invariant forms on $M$ and $M'$ which are ${\varphi}$ compatible.  They determine ${\varphi}$ compatible closed holonomy invariant currents, also denoted $\ome_T$ and $\ome'_{T'}$.  In particular, 
$$  
\lan (\xi, \xi'), (\ome_T,\ome'_{T'}) \ran \,\, = \,\,  
\lim_{s \to \infty} \left(\int_{T \ssm T_s}  \hspace{-0.5cm} \xi \wedge \ome_T - \int_{T' \ssm T_s'}  \hspace{-0.5cm}\xi' \wedge \ome'_{T'} \right).
$$
Here $\ome_T = \ome \, |_T$, which is  well defined and is holonomy invariant, as is $\ome'_{T'}$.

 For Riemannian foliations, examples of this type abound.  In particular,  the characteristic forms of holonomy invariant bundles which agree at infinity, for example $\wedge^j \nu^* \otimes (\otimes^{\ell} \nu)$, and $\wedge^j \nu'^* \otimes (\otimes^{\ell} \nu')$.  

For definiteness, we will {{generally}} use this example in the sequel, {{but all the statements obviously remain valid with more general holonomy invariant currents}}.
\end{enumerate}
\end{remark}

In this paper, we will have a number of different pairings, which will be uniformly indicated by the notation $\langle \cdot, \cdot \rangle$.  The notation should make clear where the objects live.  For example, we have 
$$
\langle\left[\int_F \AS(D_F),\int_{F'}AS(D'_{F'})\right], {{[\ome_T,\ome'_{T'}]}}\rangle \; = \; 
\int_T  \left( \int_F \AS(D_F)\right) \wedge \ome_T - \int_{T'}  \left( \int_{F'}AS(D'_{F'}) \right) \wedge \ome'_{T'},
$$
and
$$
\langle (\ch(P_0), \ch(P_0)),  (\ome_T, \ome'_{T'})  \rangle \; = \; 
\int_T \ch(P_0) \wedge \ome_T \, - \, \int_{T'} \ch(P'_0) \wedge \ome'_{T'} 
$$
In the first case, the terms in the pairing live in relative Haelfliger cohomology.  In the second, the terms are pairs of bounded Haefliger forms, and the second pair happen to agree near infinity.

\section{Chern characters in Haefliger cohomology} \label{ChHae}

We recall in this section the main steps in the construction of the Chern character in Haefliger cohomology and explain how they immediately extend to the case of a pair of foliations which are compatible near infinity. In this latter case, our Chern character takes values in a relative version of Haefliger cohomology that we introduce below.

\ms

In \cite{BH21} we worked on $M$, while in \cite{H95,HL99,BH04, BH08}, we worked on  $\cG$, which we will also do here, but our basic data will be taken from the ambient manifolds.  The results in \cite{BH21} extend readily to $\cG$ with the only change being that the spectral projections used on $\cG$ are for the operator lifted to $F_s$.  This represents another extension, in the spirit of  Connes' extensions in \cite{C79, C81}, of the classical Atiyah $L^2$ covering index theorem,  \cite{A76}.   In addition, as will be explained below,} the results in the above cited papers where $M$ was assumed to be compact still hold provided both $M$ and $F$ are of bounded geometry and we use our definition of the Haefliger cohomology.  

\ms

All the data in the previous section may be lifted to $(\cG, F_s)$ using the map $r:\cG \to M$.  The notation we will use is obtained from that above by:
$$
E_M \to E;  \quad \cS_F \to \cS; \quad  \nabla^{F,E} \to \nabla;  \quad L \to \wtL_x; \quad  D_F \to D;  \quad D_L \to D_x.
$$
Thus the smooth {{$\cG$ invariant family}} $D= \{D_x\}$ of leafwise Dirac operators acting on sections of ${\cS \otimes E} |\wtL_{x}$ is given as follows.  Let $X_1,\ldots, X_p$ be a local oriented orthonormal basis of $T\wtL_{x}$.  Then,
$$
D_x=\sum^p_{i=1}\rho(X_i)\nabla_{X_i} \; :\; C^{\infty}_{c}(\maG_x, {\cS}^{\pm} \otimes E) \to C^{\infty}_{c}(\maG_x, {\cS}^{\mp} \otimes E)\;  
\text{ and } \;
D_x^2: C^{\infty}_{c}(\maG_x, {\cS}^{\pm} \otimes E) \to C^{\infty}_{c}(\maG_x, {\cS}^{\mp} \otimes E). 
$$
Denote by $\wedge \nu_s^*$, the exterior powers of the dual normal bundle $\nu_s^*$ of $F_s = r^*F,$ which we identify with $s^{*}(T^{*}M) = s^*(TF^*) \oplus s^*(\nu^*)$ so that each $C^{\infty}_{c} (\cS \otimes E\otimes \wedge \nu_s^*)$ is an $\Omega^*(M)$-module.
We extend $D$ to an $\Omega^*(M)$-equivariant operator
$$
D: C^{\infty}_{c}(\cS \otimes E\otimes \wedge \nu_s^*)   \longrightarrow 
C^{\infty}_{c}(\cS \otimes E\otimes \wedge \nu_s^*),
$$
by using the leafwise flat connection on $\wedge \nu_s^*$ determined by the pull-back of the Levi-Civiti connection on 
$T^{*}M$. 

\ms

 In \cite{BH08}, we used the traces $\Tr$ and  $\hTr$ to define Connes-Chern characters in $H^*_c(M/F)$ for operators on $C^{\infty}_{c}({{\cS \otimes E}})$. For the leafwise spectral projection $P_0$ onto the kernel of $D^{2}$, when this latter is smooth,  this is denoted,
$$
{{\ch}} (P_{0}) \; \in  \; \cohom_c^*(M/F).
$$
We also proved that if $M$ is compact and $\Ind_{a}(D)$ is Connes' $K$-theory index class defined in terms of a parametrix for $D$, then under the usual regularity assumption, ${{\ch}} (P_{0}) =  {{\ch}} (\Ind_{a}(D))$.  We now extend these notions to our situation.

\ms

We now return to our compatible foliations $(M, F)$ and $(M', F')$ and their  holonomy groupoids $\cG$ and $\cG'$.  First, we lift the compatibility data $\Phi$ to $\cG$ and denote again the corresponding data by $\Phi$, which gives an equivalence off  (the generally non-compact subsets) $K = r^{-1}(K_M)$ and $K' = (r')^{-1}(K_{M'})$, that is on the subsets $V = r^{-1}(V_M)$ and $V' = (r')^{-1}(V_{M'})$.  In \cite{BH08}, we defined an algebra of super-exponentially decaying $\cG-$operators on $C^{\infty}_{c}(\cS \otimes E\otimes \wedge \nu_s^*)$.  Here we need a stronger condition on our operators, namely that they have finite propagation.   This is provided by using operators from the algebra $C_u^{\infty}( (\cS \otimes E\otimes \wedge \nu_s^*) \boxtimes (\cS \otimes E\otimes \wedge \nu_s^*)^*)$, which we denote simply as $C_u^{\infty}( F_s)$.  Any $A = (A_x)_{x\in M} \in    C_u^{\infty}( F_s)$ defines a leafwise (smoothing) $\cG-$operator on
$C^{\infty}_{c}(\cS \otimes E\otimes \wedge \nu_s^*)$ which has uniform finite propagation, and its Schwartz kernel  is smooth in all variables, with all derivatives being globally bounded, the bounds possibly depending on the derivatives.

\ms

Using the algebra $C_u^{\infty}( F_s)$, we have a $K$-theory index class represented by idempotents constructed from a parametrix, and this  $K$-index does not depend on the parametrix, so its Connes-Chern character is also independent of the parametrix. 
In particular, as $D$ is an odd super operator, we may write 
$
D \, = \, \left[\begin{array}{cc}  0& D^- \\ D^+  & 0\end{array}\right].
$
Suppose that $Q_t$ is a smooth (in $t$) family of  leafwise parametrices for $D$.  That is, each $Q_t$ is an odd operator which is smooth in all variables, and which has finite propagation remainders, namely the even operators 
$$
S_t \; = \; \Id_{\cS^+ \otimes E} \, - \, {Q_t }^-D^+  \quad  \text{ and } \quad R_t \; = \; \Id_{\cS^- \otimes E} - \, D^+{Q_t }^-.
$$
For $t > 0$,  set, as in \cite{BH08},
$$
A_t = \left[\begin{array}{cc} S_t^2 & Q_t ^-(R_t + R_t^2)
\\ &   \\
R_tD^+ &  -R_t^2 \end{array} \right].
$$
Then $A_t$ has finite propagation, is smooth in all variables, and is a bounded leafwise smoothing operator, that is,
$A_t \in C_u^{\infty}(F_s)$. 
Set $\pi_- =  \diag(0,\Id_{\cS^- \otimes E}), $ and $\pi_+ = \diag(\Id_{\cS^+ \otimes E}, 0)$.
Then $A_t + \pi_-$ is an idempotent as is $\pi_-$.  Set 
$$
\Ind_a (D) \; = \; [A_t + \pi_-] - [\pi_-] \in  \oK_0(C_u^{\infty}(F_s)).
$$ 
Since
$A_t + \pi_-$ is a smooth family of idempotents, it follows from results of \cite{BH04} that $\Ind_a (D)$ is
independent of $t$.  Since any two parametrices can be joined in a smooth family, it follows immediately that $\Ind_a (D)$ does not depend on the parametrix.

\ms

 For details of the following, see \cite{BH08}, Section 3, where we define the quasi-connection, 
$$
C^{\infty}(\cS \otimes E\otimes \wedge \nu_s^*)\stackrel{\nabla^{\nu}}{\longrightarrow} 
C^{\infty}(\cS \otimes E\otimes \wedge \nu_s^*).
$$
Given an operator $A$ on $\cS \otimes E\otimes \wedge \nu_s^*$, denote by 
$$
\pa_{\nu}: \End(C^{\infty}(\cS \otimes E\otimes \wedge \nu_s^*)) \to
\End(C^{\infty}(\cS \otimes E\otimes \wedge \nu_s^*))
$$ 
the linear operator given by the graded commutator
$$
\pa_{\nu} (A) = [\nabla^{\nu}, A].
$$
Set  $\theta =(\nabla^{\nu})^2$, which is a leafwise differential operator with coefficients in $\wedge \nu^*_s$.   Since $\pa_{\nu}^2$ is not necessarily zero, we {{used}} Connes' X-trick in \cite{BH08} to construct a new differential operator $\delta$ out of $\pa_{\nu}$ and $\theta$, whose square is zero.  Note carefully that $\delta A$ is nilpotent since it always contains a coefficient from $\wedge\nu_s ^{* \geq 1}$.

\ms

Corollary 3.7 of \cite{BH08} states,
\begin{prop}\label{Cor3.7} The Haefliger form $\Tr \Bigl{(}  A_t \exp \left[{\frac{-(\delta A_t)^2}{2i\pi}}\right]  \Bigr{)}
$ is closed, and the Haefliger class \\  $\hTr \Bigl{(}  A_t \exp \left[{\frac{-(\delta A_t)^2}{2i\pi}}\right]  \Bigr{)}$ is independent of $t$.
\end{prop}

\begin{defi}  The {{Connes-Chern character}} of $\Ind_a(D)$ is,
$$
{{\ch}}(\Ind_a(D)) = \hTr \Bigl{(}  A_t \exp \left[{\frac{-(\delta A_t)^2}{2i\pi}}\right]  \Bigr{)} \;  \in \; H_c^*(M/F).
$$
\end{defi}

We have the same constructions for $D'$.  In Section \ref{proofs}, we construct families of parametrices $Q_t$ and  $Q'_t$ directly from $D$ and $D'$ in such a way that their remainders are $\Phi$ compatible, so also are $A_t$ and $A_t'$. 

\ms

For pairs $(A, A')$ of operators from $C^\infty_u (F_s)\times C^\infty_u (F'_s)$ which are $\Phi$-compatible, there is also an algebra $C^\infty_u (F_s, F'_s; \Phi)$, and the previous construction of the analytic index class extends immediately to yield the relative analytic index class
$$
\Ind_a (D, D') \; = \; [(A_t + \pi_-, A'_t+\pi'_-] - [(\pi_-, \pi'_-)] \in  \oK_0(C_u^{\infty}(F_s, F'_s; \Phi)).
$$
The Connes-Chern character then extends to the relative case 
$$
\ch:  \oK_0(C_u^{\infty}(F_s, F'_s; \Phi)) \longrightarrow H^*_c (M/F, M'/F';\varphi),
$$
with the obvious definition  (see \cite{BH08}, Theorem 3.2  for the notation below and more precise details),
$$
\ch([\tilde{e}, \tilde{e}']) \; = \; \left[ \Tr \left( e \exp \left(\frac{-(\delta e)^2}{2i\pi}\right)\right), 
\Tr \left( e' \exp \left(\frac{-(\delta e')^2}{2i\pi}\right)\right) \right] \; \in \; H^*_c (M/F, M'/F';\varphi).
$$

\begin{definition}  Suppose the parametrices $Q_t$ and $Q'_t$ have $\Phi$ compatible remainders, so with $\Phi$ compatible operators $A_t$ and $A'_t$.  Then
the relative Connes-Chern character of $\Ind_a(D, D')$ is given by
$$
{{\ch}}(\Ind_a (D, D')) \; = \; \left[ \Tr \left( A_t \exp \left(\frac{-(\delta A_t)^2}{2i\pi}\right)\right), 
\Tr \left( A_t' \exp \left(\frac{-(\delta A_t')^2}{2i\pi}\right)\right) \right] \; \in \; H^*_c (M/F, M'/F';\varphi).
$$
\end{definition}
The class ${{\ch}}(\Ind_a (D, D'))$  is clearly well defined due to its independence of the $\Phi$-compatible pair of finite propagation parametrices. This is proved below, see Theorem \ref{Hformsclosed}, where we also point out that it is independent of the parameter $t$.

\section{Four Theorems}\label{RiemFols}

Our first main theorem is the following extension of a classical Atiyah-Singer Index Theorem.  This theorem is purely local and, as in \cite{BH21},  requires bounded geometry.   

\ms

Denote by $[AS(D_F)]$ the Atiyah-Singer characteristic class for $D_F$, and similarly for $D'_{F'}$. 
Note that for large $s$,  the differential forms satisfy $AS(D_F) = \varphi^*(AS(D'_{F'}))$ on $M_s$, so
$$
\left( \int_F \AS(D_F),\int_{F'}AS(D'_{F'})) \right) \; \in \; \maA_c^* (M/F, M'/F'; \varphi).
 $$
\begin{definition}   The relative topological index of  $(D, D')$ is,
$$
\Ind_t (D, D') \; = \; \left[\int_F \AS(D_F),\int_{F'}AS(D'_{F'}))\right]  \; \in \; H^*_c (M/F, M'/F';\varphi).
$$
\end{definition}

\begin{theorem}\label{classicalindex} [The Higher Relative Index Theorem]
Suppose that  $(M,F)$, $(M',F')$, $D$ and $D'$ are as in Section \ref{setup}.  In particular, $F$ and $F'$ need not be Riemannian.  Then,  
$$
\ch(\Ind_a (D, D')) = \Ind_t (D, D')\;  \in H^*_c \; (M/F, M'/F';\varphi)
$$
In particular, for any closed $\varphi$-compatible pair $(C, C')$ of holonomy invariant closed Haefliger currents, the following scalar formula holds
$$
\langle {{\ch}}(\Ind_a (D, D')), [C, C']\rangle = \lim_{s\to +\infty} \left(\langle \int_F \AS(D_F)\vert_{T\smallsetminus T_s}, C\rangle - \langle \int_{F'} \AS(D'_{F'})\vert_{T'\smallsetminus T'_s}, C'\rangle\right).
$$.
 \end{theorem}

Denote by $P_{(0,\ep)}$ the spectral projection for $D^{2}$ for the interval $(0,\ep)$. 
The Novikov-Shubin invariants $NS(D)$ of
$D$ are greater than $k\geq 0$ provided that there is $\tau > k$ so
that
$$
\Tr(P_{(0,\ep)}) \mbox{  is  
}  {\cO}(\epsilon^{\tau})   \mbox{  as  }   \epsilon \to 0.
$$
A Haefliger form $\Psi$ depending on $\epsilon$ is
${{\cO}}(\epsilon^\tau)$ as $\epsilon \to 0$ means that there is a
representative $\psi \in \Psi$ defined on a transversal $T$, and a constant $C > 0$, 
so that the function on $T$, $\|\psi\|_T \leq
C\epsilon^\tau$ as $\epsilon \to 0$.  Here $\| \quad \|_T$ is the
pointwise norm on forms on the transversal $T$ induced from the
metric on $M$.  

\ms

Recall that $P_0$ is the spectral projection onto the kernel of $D^{2}$.  In general the leafwise operators $P_{(0,\ep)}$ and $P_0$ are not transversely smooth (although they are always leafwise smooth), so that, in general, their Haefliger traces in $\cA^*_c(M/F)$ {{are not defined}}.  When $P_0$ is transversely smooth, see \cite{BH08}, Definition 3.8, 
$$
{{\ch}}(P_0)  \; = \; \hTr \Bigl{(}  \pi_{\pm}P_0 \exp \left({\frac{-(\delta (\pi_{\pm} P_0))^2}{2i\pi}}\right)  \Bigr{)} \; \in \; H_c^*(M/F),
$$
and similarly for $P_0'$.   Here $\pi_{\pm}$ is the grading operator 
$$
\pi_{\pm} \, = \, \diag(\Id_{\cS^+ \otimes E}, - \Id_{\cS^- \otimes E}).
$$
When $P_{(0,\ep)}$ is transversely smooth,  
$$
{{\ch}}(P_{(0,\ep)})  \; = \; \hTr \Bigl{(}  \pi_{\pm}P_{(0,\ep)} \exp \left({\frac{-(\delta (\pi_{\pm} P_{(0,\ep)}))^2}{2i\pi}}\right)  \Bigr{)} \; \in \; H_c^*(M/F),
$$  
and similarly for $P'_{(0,\ep)}$
For simplicity of notation, we will uniformly suppress the constant $2i\pi$ in what follows.   As the closed Haefliger differential forms  $\Tr(\pi_{\pm}P_0 \exp \left({-(\delta (\pi_{\pm} P_0))^2}\right))$ and $\Tr(\pi'_{\pm}P'_0 \exp \left({-(\delta (\pi'_{\pm} P'_0))^2}\right))$ are not ${\varphi}$ compatible in general, we proceed as follows.  

\ms

The component of ${{\ch}}(\Ind_a(D, D'))$ in $H^{2k}_c (M/F, M'/F';\varphi)$ is denoted $\ch^k(\Ind_a(D, D'))$, and the part of  ${{\ch}}(\Ind_a(P_0)$ in $H^{2k}_c (M/F)$ is denoted $\ch^k(P_{0}))$, and similarly for $P'_0$.

\ms

The  following theorem generalizes the main result of \cite{BH08} to bounded geometry foliations.
\begin{theorem}\label{Secondhalf1}[Riemannian Foliation Relative Index Bundle Theorem]
Fix $0 \leq \ell \leq q/2$, where $q$ is  the codimension of $F$ and $F'$.  Assume that:
\begin{itemize} 
\item the foliations $F$ and $F'$ are Riemannian;
\item the leafwise operators $P_0$, $P_0'$, $P_{(0,\ep)}$  and $P_{(0,\ep)}'$ (for $\epsilon$ sufficiently small) 
are transversely smooth;
\item  $NS(D)$ and $NS(D')$ are greater than $\ell$.
\end{itemize}
Then, for $0 \leq k \leq \ell$,  we have in $H_c^{2k}(M/F)\times H_c^{2k}(M'/F')$
$$
(\pi\times \pi') \ch^k(\Ind_a(D, D'))\, = \,   (\ch^k(\Ind_a(D)), \ch^k(\Ind_a(D'))) \, = \, 
(\ch^k(P_{0}), \ch^k(P'_{0})).
$$
\end{theorem}

\begin{remarks}\label{remarks4.4}\
\begin{enumerate} 
\item If the foliations $F$ and $F'$ are not Riemannian then we can still prove this equality but under  the stronger assumption  that $NS(D)$ and $NS(D')$ be greater than $3q$, see \cite{HL99, BH23}.
\item The examples in \cite{BHW14} show that the conditions on the Novikov-Shubin invariants are the best possible.
\item Note that if there are uniform gaps in the spectrums at $0$, that is there is $\ep >0$ so $P_{(0,\ep)}= P_{(0,\ep)}' =0$, then $P_{(0,\ep)}$  and $P_{(0,\ep)}'$ are transversely smooth and the Novikov-Shubin invariants are infinite.   For top dimensional foliations, i.e. $TF = TM$, these special cases were studied for instance in \cite{Vi67, Do87}. 
\end{enumerate}  
\end{remarks}
 
Combining Theorem \ref{classicalindex} and Theorem \ref{Secondhalf1}, we immediately  deduce the following important corollary.

\begin{theorem}\label{HaefTheorem}
Under the assumptions of Theorem \ref{Secondhalf1}, assume furthermore that $P_0=P'_0=0$, then
$$  
\left(\int_F \AS(D_F),\int_{F'}AS(D'_{F'}) \right)  \, = \,  (0,0) \quad \text{in $H_c^*(M/F) \times H_c^*(M'/F')$.}
$$ 
\end{theorem}
So the vanishing conclusion of the previous theorem holds in particular when there exists $\ep >0$ such that $P_{[0, \ep)}=0$ and $P'_{[0, \ep)}=0$.

\ms

Denote by  $\ome  \in C^{\infty}(\wedge^* \nu^*)$ and $\ome' \in C^{\infty}(\wedge^* {\nu'}^*)$ closed bounded holonomy invariant differential forms on $M$ and $M'$ which are ${{\varphi}}$ compatible.   For simplicity, we will assume that $\ome$ and $\ome'$ are $\varphi$ compatible on $V_M$ and $V'_{M'}$.  These determine $\varphi$ compatible closed bounded Haefliger forms on $T$, denoted $\ome_T$ and $\ome'_{T'}$.  Recall that $dx$ is the global volume form on $M$. 

\begin{theorem}\label{Secondhalf2}[Higher Relative Index Pairing Theorem] {{In addition to the assumptions in Theorem \ref{Secondhalf1}, assume the following:
\begin{itemize} 
\item  for $\epsilon$ sufficiently small, $P_{[0,\ep)}$ satisfies
$\dd \int_M \tr (P_{[0,\ep)}) dx \; < \; \infty$, and similarly for $P'_{[0,\ep)}$;
\item $M$, and so also $M'$, has sub-exponential growth.
\end{itemize}
Then, for any homogeneous $\ome  \in C^{\infty}(\wedge^{q-2k} \nu^*)$ and $\ome' \in C^{\infty}(\wedge^{q-2k}  {\nu'}^*)$ as above, ($0 \leq k \leq \ell$), 
$$
\dd \int_T {{\ch}}(P_0) \wedge \ome_T\text{ and }\dd \int_{T'} \ch(P'_0) \wedge \ome'_{T'}\text{  are well defined complex numbers,}
$$
and
$$
\int_T \ch(P_0) \wedge \ome_T \, - \, \int_{T'} \ch(P'_0) \wedge \ome'_{T'}  \, = \,   \langle\left[\int_F \AS(D_F),\int_{F'}AS(D'_{F'})\right],{{\left[\ome_T,\ome'_{T'}\right]}}  \rangle. 
$$}}
 \end{theorem}

\begin{remarks}\
\begin{enumerate}
\item {{Since the pair of Connes-Chern characters of $P_0$ and $P'_0$ is {{usually}} not $\varphi$-compatible,  the previous theorem is totally new and we cannot deduce it from any absolute version of the index bundle theorem. This is compatible with the classical relative index theorem.}}
\item The theorem also holds for appropriate closed $\varphi$ compatible closed holonomy invariant currents, but this more general statement will not be needed for our applications.
\item  We shall see in Section \ref{TIMAFTS} that the finite integral assumptions are satisfied when the zero-th order operator $\maR^E_F$ defined there in the Bochner formula is strictly positive near infinity.  As $\maR^E_F$ is locally defined, this means that $\maR^{E'}_{F'}$ is also strictly positive near infinity. 
\item {{The growth condition is a technical assumption which simplifies the proof, it   can  be weakened as explained in Remark \ref{growth}}}
\item {{The main theorem in \cite{BH21} recovers the Gromov-Lawson relative index theorem in full generality for bounded geometry manifolds, which correspond to top-dimensional foliations. 
Our  results here require more conditions to deal with the higher components of the Connes-Chern character, and it only recovers the Gromov-Lawson results for sub-exponential bounded geometry manifolds. Recall that in the top-dimensional case,  Gromov-Lawson show that there is $\ep > 0$ so that $P_{(0,\ep)} = 0$, and  $\dd \int_M \tr(P_0) dx < \infty$, so  all the other assumptions of Theorems \ref{Secondhalf1} and \ref{Secondhalf2} are fulfilled.}}
\end{enumerate}
\end{remarks}

\section{Proofs of the Theorems}\label{proofs}
    
This section is devoted to the proofs of Theorems \ref{classicalindex}, \ref{Secondhalf1}, \ref{HaefTheorem} and \ref{Secondhalf2}. The proofs are rather technical and have been split into many intermediate lemmas and propositions. We shall first prove Theorem \ref{classicalindex} and then later on Theorems \ref{Secondhalf1} and \ref{HaefTheorem}, and eventually we shall end this section by the proof of Theorem \ref{Secondhalf2}.

\ms

Recall the following construction from \cite{BH21}.  Denote the Fourier Transform of a complex valued function $g$ by $\what{g}$ and $FT(g)$, and its inverse transform $FT^{-1}(g)$ by $\wtit{g}$.  If $h$ is also a complex function, denote  the convolution of $g$ and $h$ by $g \star h$.  Set $g_\lambda(z) = g(\lambda z)$, for non-zero $\lambda \in \R^*$.   We have the following facts:  
$$
FT (g_\lambda) = \frac{1}{\lambda} FT (g)_{ \frac{1}{\lambda}}; \;  FT (g \star h)= \sqrt{2\pi}FT(g) FT(h); \text{  and  } FT(\what{g}) = FT^{-1}(\what{g}) = g, \text{  if $g$ is even}.
$$
Fix  a smooth even non-negative function $\psi$ supported in $[-1, 1]$, which equals $1$ on $[-1/4, 1/4]$, is non-increasing on $\R_+$, and whose integral over $\R$ is $1$.  Note that $FT( \what{\psi}) = \psi$ since $\psi$ is even.  The family $\frac{1}{\sqrt{t}} \what{\psi}_{\frac{1}{\sqrt{t}}}$ is an approximate identity when acting on a Schwartz function $f$ by convolution, since, up to the constant $\sqrt{2\pi}$ which we systematically ignore,
$$
\frac{1}{\sqrt{t}} \what{\psi}_{\frac{1}{\sqrt{t}}} \star f  \; = \; FT^{-1}(FT( \frac{1}{\sqrt{t}} \what{\psi}_{\frac{1}{\sqrt{t}}} \star f ))  \; = \; FT^{-1} (\psi_{\sqrt{t}}  \what{f}) \to \wtit{\what{f}} \; = \; f,
$$
in the Schwartz topology as $t \to 0$.  Denote as usual by $|| \,\cdot\, ||_{r,s}$  the norm of an operator acting from the $r$ Sobolev space to the $s$ Sobolev space.  Then more is true.

\begin{lemma}\label{Schwlimit}
Suppose that $\mu:\R_+\to \R_+$,  with $\mu (t)\leq C_p t^p$ or $\mu (t)\geq C_p t^{-p}$ near $0$,  where $p > 0$ and $C_p >0$.  Then, for any Schwartz function $f$,
$$
 \lim_{t\to 0} \left(\left[\frac{1}{\sqrt{t}} \what{\psi}_{\frac{1}{\sqrt{t}}} \star f\right]_{\mu (t)} \; -  \;   f_{\mu (t)}\right)   =  0
$$
 in the Schwartz topology.
 
Thus for all $r,s$, 
$$
\lim_{t\to 0} ||\left[\frac{1}{\sqrt{t}} \what{\psi}_{\frac{1}{\sqrt{t}}} \star f\right]_{\mu (t)}(D) \; -  \;   f_{\mu (t)}(D)||_{r,s} \; =  \; 0,
$$
so the differences of their Schwartz kernels converge uniformly to  $0$ pointwise. 
\end{lemma}
 
\begin{proof} The last statement follows from standard Sobolev theory given the first.
Thus we need only prove that the difference of the Fourier transforms goes to zero in the Schwartz topology. But, 
$$
FT \left(\left[\frac{1}{\sqrt{t}} \what{\psi}_{\frac{1}{\sqrt{t}}} \star f\right]_{\mu (t)}\right) - FT \left[f_{\mu (t)}\right]  =  
\frac{1}{\mu (t)} \what{f}_{\frac{1}{\mu (t)} }   (\psi_{\sqrt{t}/\mu (t)} -1).
$$
Now, $\psi_{\sqrt{t}/\mu (t)}(z)-1$ is $0$ for $|z| \leq\mu (t)/4\sqrt{t}$ and constant for $|z| \geq \mu (t)/\sqrt{t}$, so all its derivatives are zero on these subsets.  In addition, for all non-negative $n$, there is a constant $C_n$ so that 
$$
|\frac{\pa^n}{\pa z^n} (\psi_{\sqrt{t}/\mu (t)}(z)-1)| \; \leq \; C_n (\sqrt{t}/\mu (t))^n.
$$
 Thus, we have  
$$
\vert\vert z^n \frac{\pa^m}{\pa z^m} \left[ \frac{1}{\mu (t)} \what{f}_{\frac{1}{\mu (t)} } (\psi_{\sqrt{t}/\mu (t)} -1)\right] \vert\vert_\infty \; = \; 
\sup_{\vert z\vert \geq \mu (t)/(4\sqrt{t})} \left\vert z^n \frac{\pa^m}{\pa z^m} \left[ \frac{1}{\mu (t)} \what{f}_{\frac{1}{\mu (t)} } (\psi_{\sqrt{t}/\mu (t)} -1)\right] \right\vert \;  \leq  \;
$$
$$
 \sup_{\vert z\vert \geq \mu (t)/(4\sqrt{t})}\left\vert z^n \sum_{k=0}^m C_{m-k} (\sqrt{t}/\mu (t))^{m-k} \frac{\pa^k}{\pa z^k} \left[\frac{1}{\mu (t)} \what{f}_{\frac{1}{\mu (t)} }  \right] (z) \right\vert 
\; =  \;   
$$
$$\sum_{k=0}^m C_{m-k} (\sqrt{t}/\mu (t))^{m-k}
 \mu (t)^{-(k+1)} \sup_{\vert z\vert \geq  \mu (t)/(4\sqrt{t})} \left\vert z^n \what{f}^{(k)} \left(z /\mu (t)\right) \right\vert \\
\; = \;
$$
$$ 
\sum_{k=0}^m C_{m-k}(\sqrt{t})^{m-k}\mu (t)^{-(m+1)}  \sup_{\vert z\vert \geq 1/(4\sqrt{t})} \left\vert z^n \what{f}^{(k)} (z) \right\vert. 
$$

Since $f$, so also $\what{f}$, is Schwartz, for any non-negative $k\in \Z$, the function $z\mapsto z^n \what{f}^{(k)} (z)$ is Schwartz.   But for any Schwartz function $g$, any $N\geq 0$ ($N\leq 0$ is trivial) and any $\eta >0$, $\dd
\lim_{t\to 0} t^{-N} \hspace{-0.3cm}\sup_{\vert z\vert \geq \eta/\sqrt{t}} \vert g(z)\vert =0.$ Thus, if $\mu (t)\leq C_p t^p$ or $\mu (t)\geq C_p t^{-p}$ near $0$,  

\noi \hspace{4.5cm} 
$\dd\lim_{t \to 0}\vert\vert z^n \frac{\pa^m}{\pa z^m} \left[ \frac{1}{\mu (t)} \what{f}_{\frac{1}{\mu (t)} }  (\psi_{\sqrt{t}/\mu (t)} -1)\right] \vert\vert_\infty \; = \; 0.$
\end{proof}

Define the functions $\alp(t)$ and $\beta(t)$ as follows.  Both have domains  $(0,1)$, and are smooth.  $\alp(t) = t$ near $0$, and $\alp(t) = 1-t$ near $1$, it is increasing on $(0,1/2]$ and symmetric about $t = 1/2$.   $\beta$ is an increasing function, with $\beta = t$ near $0$, and $\beta(t) = (1-t)^{-1}$ near $1$. 

\ms

Set  $e(z) = e^{-z^2/2}$, and for $t \in (0,1)$, set  
$$
\chi^t(z) \; =  \; \left[\frac{1}{\sqrt{\alp(t)}} \what{\psi}_{{\frac{1}{\sqrt{\alp(t)}}}}\star e  \right]_{\sqrt{\beta(t)}} (z).
$$
\begin{remark}\label{limseqs}
By Lemma \ref{Schwlimit}, we have, 
$$
\lim_{t \to 0}\left(\chi^t(z) \; - \; e^{-tz^2/2}\right) \; = \;  0 \; =  \;\lim_{t \to 1}\left(\chi^t(z) \; - \; e^{-z^2/(2(1-t))}\right),
$$  
in the Schwartz topology.  In addition, the limit as $t \to 0$ of the Schwartz kernel of $\chi^t(D) \; - \; e^{-tD^2/2}$ and  the limit as $t \to 1$ of the Schwartz kernel of $\chi^t(D) \; - \; e^{-D^2/(2(1-t))}$ both converge uniformly pointwise to zero.
\end{remark}

\begin{lemma}\label{finprop}
$\chi^t(D)$ has finite propagation $ \leq \sqrt{\beta(t)/\alp(t)}$. 
\end{lemma} 
\begin{proof}
Since $\what{e} = e$, we have that up to a constant,
$$
FT(\frac{1}{\sqrt{\alp(t)}}  \what{\psi}_\frac{1}{\sqrt{\alp(t)}}  \star e)  \; = \;  \psi_{\sqrt{\alp(t)}}  e.  
$$
In fact, up to a constant,
$$
\chi^t (D)   \; = \;  FT^{-1}(\psi_{\sqrt{\alp(t)}}  e)(\sqrt{\beta(t)}D)   \; = \; 
\int_\R \psi (\sqrt{\alp(t)} \xi)  e (\xi) \cos (\xi \sqrt{\beta(t)} D) \, d\xi,
$$
since $\psi_{\sqrt{\alp(t)}} e$ is even.  Setting $\eta = \sqrt{\alp(t)}\xi $ 
gives,
$$
\chi^t (D) = \frac{1}{\sqrt{\alp(t)}} \int_{\vert \eta\vert \leq 1} \psi (\eta)  e (\eta/\sqrt{\alp(t)}) \cos (\eta\sqrt{\beta(t)/\alp(t)} D)  d\eta.
$$
The operator $\cos (\eta \sqrt{\beta(t)/\alp(t)}D)$ has propagation $\leq \vert \eta  \sqrt{\beta(t)/\alp(t)}\vert$, see \cite{Ch73, Roe87}.  Thus $\chi^t (D)$ has finite propagation $\leq \sqrt{\beta(t)/\alp(t)}$, which near $0$ is $\leq 1$, while near $1$ it is $\leq (1-t)^{-1}$, so may go to infinity as $t \to 1$.
\end{proof}

As $D$ is an odd super operator, we may write 
$$
D \, = \, \left[\begin{array}{cc}  0& D^- \\ D^+  & 0\end{array}\right],  \text{ and we set } 
Q_t  \, = \,    \Big(\frac{1- \chi^t (0)^{-1}\chi^t (z)}{ z}\Big)(D) \, = \,    
\Big(\frac{1- \chi^t (0)^{-1}\chi^t (z)}{ z^2}z\Big)( D).
$$
We claim that $Q_t$ is a smooth family of  leafwise parametrices for $D$ with finite propagation $\Phi$ compatible remainders, namely the even operators 
$$
S_t \; = \; \Id_{\cS^+ \otimes E} \, - \, {Q_t }^-D^+  \quad  \text{ and } \quad R_t \; = \; \Id_{\cS^- \otimes E} - \, D^+ Q_t ^-.
$$
There are similar relations for $D^-$.  

\ms

The main step in the proof of Theorem \ref{classicalindex} is the following expected independent result.
\begin{prop}\label{Ptdefn} For $0 < t  < 1$,  set, as in \cite{BH21},
$$
A_t = \left[\begin{array}{cc} S_t^2 & Q_t ^-(R_t + R_t^2)
\\ &   \\
R_tD^+ &  -R_t^2 \end{array} \right] \; =
\left[\begin{array}{cc} S_t^2 & S_tQ_t ^-(1 + R_t)
\\ &   \\
R_tD^+ &  -R_t^2 \end{array} \right],
$$
a form which is more useful here.
Then $A_t$ and $\delta A_t$, so also $(\delta A_t)^2$,  have finite propagations which are bounded by multiples of $\sqrt{\beta(t)/\alp(t)}$,  are smooth in all variables, and are bounded leafwise smoothing operators. 
\end{prop}  

\ms

\begin{proof} 
We deal with $A_t$ first.  
Note that $S_t = \chi^t (0)^{-1} \chi^t(D)$ acting on $\cS^+ \otimes E$, and similarly for $R_t$ acting on $\cS^- \otimes E$.  They both have finite propagations, and by Theorem 2.1, \cite{Roe87},  they are both smooth in all variables.  It follows immediately that $S_t^2$, $R_t^2$, $R_tD^+$, $S_t$ and $R_t$ are also smooth in all variables.  Since propagation is additive for compositions, they all  have finite propagations, which are bounded by multiples of $\sqrt{\beta(t)/\alp(t)}$.  Finally, since $\chi^t(z)$ is a Schwartz function, $\chi^t(D)$ and $\chi^t(D)D$ are bounded leafwise smoothing operators.

\ms

To deal with $Q_t ^-$, we show that
$$ 
\wtit{Q}_t(D^-D^+) \; = \;   \frac{\left[\left(1-  \chi^t (0)^{-1}\chi^t (z)\right)( D)\right]^+}{ D^-D^+} 
$$
has finite propagation which is bounded by a multiple of $\sqrt{\beta(t)/\alp(t)}$, so also does $Q_t ^- = \wtit{Q}_t(D^-D^+) D^-$, and that $S_tQ_t^-$ is smooth in all variables and is a bounded leafwise smoothing operator.

\ms

For $u \in (0,1]$, set   
$$
\chi^{t,u}(z) \, =  \, \left[\frac{1}{\sqrt{\alp(t)}} \what{\psi}_{{\frac{1}{\sqrt{\alp(t)}}}}\star e_u  \right]_{\sqrt{\beta(t)}} (z),
$$  
and 
$$
\wtit{q}^{t,u} (z)    \, = \,      \frac{1- \chi^{t,u}(0)^{-1} \chi^{t,u}(z)}{z^2}
\, = \,  \chi^{t,u}(0)^{-1}    \frac{\chi^{t,u}(0) -\chi^{t,u}(z)}{z^2}.
$$
Notice that $\chi^{t,0}(z)$ and $\wtit{q}^{t,0}(z)$ are also well defined, and that for fixed $z$, the resulting function  is continuous on $[0, 1]$ and smooth on $(0, 1)$. Since  
$$
\chi^{t,u}(z) \, =  \, FT^{-1}(FT\left(\frac{1}{\sqrt{\alp(t)}} \what{\psi}_{{\frac{1}{\sqrt{\alp(t)}}}}\star e_u  \right))(\sqrt{\beta(t)}z), 
$$
we have, 
$$
\chi^{t,u}(z) \, =  \,   
\int_{\R}  \psi(\sqrt{\alp(t)}y) \frac{1}{u} e^{-y^2/2u^2} \cos(y\sqrt{\beta(t)}z) dy  \, =  \,  
\int_{\R}  \psi(\sqrt{\alp(t)}uy)  e^{-y^2/2} \cos(uy\sqrt{\beta(t)}z) dy,
$$
and,
$$
\chi^{t,u}(0)  \; = \;  \int_{\R} \psi(\sqrt{\alp(t)}uy)  e^{-y^2/2} dy.
$$
The latter is smooth in $t$ and $u$, positive, and bounded by $\dd \int_{\R}  e^{-y^2/2} dy$.  Thus $\chi^{t,u}(0)^{-1}$ is smooth on $(0,1) \times (0,1)$.
In addition,
$$  
-\frac{\pa}{\pa u}\left(\chi^{t,u}(0)^{-1} \right)  \; = \;  
\chi^{t,u}(0)^{-2} \int_{\R} \psi'(\sqrt{\alp(t)}uy)  \sqrt{\alp(t)}ye^{-y^2/2} dy,
$$
so has the same properties as $\chi^{t,u}(0)^{-1}$.

\ms

Next, we have 
$$
\wtit{q}^{t,u} (z)    \, = \,    \chi^{t,u}(0)^{-1}   \int_{\R} \psi(\sqrt{\alp(t)}uy)e^{-y^2/2 }\frac{ (1 - \cos(uy\sqrt{\beta(t)}z))}{z^2}dy.    
$$
So,
$$
\frac{\pa \wtit{q}^{t,u}}{\pa u}(z)   \, = \,    \frac{\pa}{\pa u}\left(\chi^{t,u}(0)^{-1} \right)  \int_{\R} \psi(\sqrt{\alp(t)}uy)e^{-y^2/2 }\frac{ (1 - \cos(uy\sqrt{\beta(t)}z))}{z^2}dy \; + \; 
$$
$$
 \chi^{t,u}(0)^{-1}  \int_{\R}  \sqrt{\alp(t)}\psi'(\sqrt{\alp(t)}uy) ye^{-y^2/2 }\frac{ (1 - \cos(uy\sqrt{\beta(t)}z))}{z^2}dy \; + \; 
$$
$$
 \chi^{t,u}(0)^{-1} \int_{\R} \sqrt{\beta(t)}\psi(\sqrt{\alp(t)}uy)ye^{-y^2/2}\frac{sin(uy\sqrt{\beta(t)}z)}{z}dy.
$$
For $t$ and $z$ fixed, this is a smooth function of $u$.

\ms

Note that $\chi^{t,0}(z)$ is well defined and equals $\dd \int_{\R}e^{-y^2/2} dy$, which is independent of $z$.  Thus,  $\wtit{q}^{t,0} (z)$ is also  well defined and equals $0$.  As $\wtit{q}^{t,1} (D)^+ = \wtit{Q}_t (D^-D^+)$, we have
$$
\wtit{Q}_t (D^-D^+) = \wtit{q}^{t,1} (D)^+ - \wtit{q}^{t,0}(D)^+ =  \left[\int_0^1 \frac{\pa \wtit{q}^{t,u}}{\pa u} (D) \, du\right]^+,
 $$
 so we need to show that $\dd \frac{\pa \wtit{q}^{t,u}}{\pa u}(D)$ has finite propagation which is  bounded by a multiple of $\sqrt{\beta(t)/\alp(t)}$.  Since $\chi^{t,u}(0) $ and $\frac{\pa}{\pa u}\left(\chi^{t,u}(0)^{-1} \right)$ are independent of $z$, they give multiples of the identity map when evaluated at $D$, so have propagation zero and may be disregarded.  Thus, we may assume that $\dd \wtit{q}^{t,u} (z) \, = \,     \frac{\chi^{t,u}(z) - \chi^{t,u}(0)}{z^2}$.

\ms

Since $\chi^{t,u}(z)$ is an even function, it has Taylor expansion in $z$ with integral remainder
$$
\chi^{t,u}(z) \; = \; \chi^{t,u}(0)  \, + \,     \frac{(\chi^{t,u})^{(2)}(0)}{2} z^2   \, + \,    \frac{z^4}{6} \int_0^1 (1-v)^3 (\chi^{t,u})^{(4)} (vz) dv.
$$
So the Taylor expansion in $z$ with integral remainder of $\dd \wtit{q}^{t,u} (z)  \, = \,  \frac{\chi^{t,u}(z) - \chi^{t,u}(0)}{z^2}$ is 
$$
\wtit{q}^{t,u} (z) \, = \,  \frac{(\chi^{t,u})^{(2)}(0)}{2} + \frac{z^2}{6} \int_0^1 (1-v)^3 (\chi^{t,u})^{(4)} (vz) dv.
$$
The term $\frac{\pa}{ \pa u} ((\chi^{t,u})^{(2)}(0))$ is independent of $z$, so, as above, it may be disregarded.
Using the fact that 
$$
\chi^{t,u}(vz) \, =  \,   \int_{\R}  \psi(\sqrt{\alp(t)}y) \frac{1}{u} e^{-y^2/2u^2} \cos(y\sqrt{\beta(t)}vz) dy,
$$
we have
$$
(\chi^{t,u})^{(4)}(vz) \, =  \,   
\int_{\R}  \psi(\sqrt{\alp(t)}y) \beta(t)^2v^4y^4 \frac{1}{u}e^{-y^2/2u^2} \cos(y\sqrt{\beta(t)}vz) dy \, =  \, 
$$
$$
\left[\frac{1}{\sqrt{\alp(t)}} \what{\psi}_{{\frac{1}{\sqrt{\alp(t)}}}}\star \beta(t)^2 v^4e_u^{(4)}\right](\sqrt{\beta(t)}vz) \, = \, 
\left[\frac{1}{\sqrt{\alp(t)}} \what{\psi}_{{\frac{1}{\sqrt{\alp(t)}}}}\star \beta(t)^2v^4p e_u\right](\sqrt{\beta(t)}vz),
$$ 
where $p$ is a finite polynomial in $u$ and $z$, since $e_u(z) = e^{-u^2z^2/2}$.  \underline{Note carefully} that  
$
\frac{\partial}{\partial u}(\chi^{t,u})^{(4)}(vz)
$ 
has the same form.  

\ms

As $pe_u$ is a Schwartz function for non-zero $u$, so is $(\chi^{t,u})^{(4)}(vD)$, and
the now usual argument shows that it has  propagation $\leq v\sqrt{\beta(t)/\alpha (t)}$.  As $D^2$ has zero propagation,  $\dd \frac{\pa \wtit{q}^{t,u}}{\pa u}(D)$ has  propagation which is a multiple of $\sqrt{\beta(t)/\alpha (t)}$, so also do $\wtit{Q}_t$ and $Q_t^-$.

\ms

For smoothness and bounded leafwise smoothing of $S_tQ_t ^-$,  first note that 
$$
S_t  \int_0^1 \frac{\pa}{ \pa u} \left(\frac{(\chi^{t,u})^{(2)}(0)}{2}\right)(D) du  \; = \; 
\left(  \int_0^1 \frac{\pa}{ \pa u} \left(\frac{(\chi^{t,u})^{(2)}(0)}{2}\right)(D) du\right)  S_t, 
$$
and $S_t$ has these properties.  Finally, any positive power of $D$ times an operator of the form 
$$
\left[ \frac{1}{\sqrt{\alpha(t)}} \what\psi_{\frac{1}{\sqrt{\alpha(t)}}} \star \beta(t)^2v^4p e_u\right](\sqrt{\beta (t)}vD)
$$
has these properties, (the function in the brackets is Schwartz), see for instance Theorem 2.1, \cite{Roe87}.  Thus, 
$$
\int_0^1 \frac{D^2}{6} \int_0^1 (1-v)^3 \frac{\partial}{\partial u}(\chi^{t,u})^{(4)}(vD) dv \, du
$$
has these properties.  Therefore, $S_tQ_t ^-$ has all the requisite properties, so $A_t$ does also.  

\ms

The operator $\delta A_t$ is essentially a polynomial in $A_t$, $\pa_{\nu}(A_t) = [\nabla^{\nu},A_t]$, and $\theta = (\nabla^{\nu})^2$.  Both $\nabla^{\nu}$ and $\theta$ are smooth and bounded in all variables and are differential operators. Since $A_t$ has finite propagation and is smooth in all variables,  $\delta A_t$ and  $(\delta A_t)^2$ also have finite propagations and are smooth in all variables.  

\ms

It remains to show that $\delta A_t$ is bounded leafwise smoothing, but this is a routine exercise.  We give some details for the convenience of the reader.  Every term of $\delta A_t$ contains either $A_t$, $\pa_{\nu}(A_t)$, or both.  As $A_t$ is bounded leafwise smoothing, we need only show that $\pa_{\nu}(A_t) = [\nabla^{\nu},A_t]$ is  bounded leafwise smoothing, since $\theta$ composed with a bounded leafwise smoothing operator is bounded leafwise smoothing.  As $\pa_{\nu}$ is a derivation, we need only show that $\pa_{\nu}$ applied to the individual elements of $A_t$, save $D^+$, yields a bounded leafwise smoothing operator.  

\ms

First,
$$
\pa_{\nu}(\chi^t(D)) \; =  \; 
\pa_{\nu}\left(\int_{\R} \frac{1}{\sqrt{\alp(t)}} \what{\psi}_{{\frac{1}{\sqrt{\alp(t)}}}}(y) e^{-(\sqrt{\beta(t)}D - y)^2/2} dy\right)   \; =  \; 
$$
$$
-\frac{1}{2}\int_{\R} \frac{1}{\sqrt{\alp(t)}} \what{\psi}_{{\frac{1}{\sqrt{\alp(t)}}}}(y) 
\int_0^1 e^{-(1-w)(\sqrt{\beta(t)}D - y)^2/2} \pa_{\nu}((\sqrt{\beta(t)}D - y))^2) 
e^{-w(\sqrt{\beta(t)}D - y)^2/2}dw \, dy.
$$
For the second equality, we refer to the proof of Proposition 3.5 of \cite{H95}, which is an extension of Proposition 2.8 of \cite{B86} to foliations with Hausdorff holonomy groupoids.  Now,  $\pa_{\nu}((\sqrt{\beta(t)}D - y)^2)$ is a differential operator with smooth bounded coefficients, so $\pa_{\nu}(\chi^t(D))$ has the same properties as $\chi^t(D)$, i.e.\ it is  bounded and leafwise smoothing.  Thus $\pa_{\nu}(S_t)$ and $\pa_{\nu}(R_t)$ are bounded leafwise smoothing.
Since $\pa_{\nu}(D)$ is a differential operator with smooth bounded coefficients, $R_t \pa_{\nu}(D)$ is also bounded leafwise smoothing. 
Finally, as $Q_t^-  =   \wtit{Q}_t(D^-D^+) D^-$, it suffices to show that $S_t\pa_{\nu}(\wtit{Q}_t(D^-D^+))$ is bounded leafwise smoothing.  As above, this follows if we show that  $S_t\pa_{\nu}(\frac{\pa \wtit{q}^{t,u}}{\pa u} (D))$ is bounded leafwise smoothing.  For the terms
$$
S_t\pa_{\nu}\left(\frac{\pa}{ \pa u} ((\chi^{t,u})^{(2)}(0))(D)\right) \; = \; 
\pa_{\nu}\left(\frac{\pa}{ \pa u} ((\chi^{t,u})^{(2)}(0))\right) S_t, \; \text{ and }  \; S_t \pa_{\nu}( \frac{D^2}{6})
$$ 
this is obvious.  As noted above, the term
$
\dd \pa_{\nu}\left( \int_0^1 (1-v)^3 \frac{\partial}{\partial u}(\chi^{t,u})^{(4)}(vD) dv \right)
$
has the form
$$
\pa_{\nu}\left( \int_0^1 (1-v)^3 \left[\frac{1}{\sqrt{\alp(t)}} \what{\psi}_{{\frac{1}{\sqrt{\alp(t)}}}}\star \beta(t)^2v^4p e_u\right](\sqrt{\beta(t)}vD) dv \right) \; = \; 
$$
$$
\pa_{\nu}\left( \int_0^1 (1-v)^3
\int_R \Big( \frac{1}{\sqrt{\alp(t)}} \what{\psi}_{{\frac{1}{\sqrt{\alp(t)}}}} (y) \Big) \beta(t)^2v^4p(u,\sqrt{\beta(t)}vD -y) 
e^{-u (\sqrt{\beta(t)}vD - y)^2/2} dy \,dv \right).
$$
The argument used for $\pa_{\nu}(\chi^t(D))$ is also valid here, so we have the result.
\end{proof}
We have the same results for $D'$, and since $A_t$ and $A_t'$ are constructed directly from $D$ and $D'$ {{and have finite propagation,}} they are $\Phi$ compatible, as are $\delta A_t$ and $\delta A_t'$.  Thus $\Tr \left( A_t\exp(-(\delta A_t)^2\right)$ and $\Tr \left(  A'_t\exp (-(\delta A'_t)^2 \right)$ are $\varphi$ compatible.  Now Theorem \ref{classicalindex} will be deduced right away from the following

\begin{theorem} \label{Hformsclosed} 
For $t \in (0,1)$, the  $\varphi$ compatible Haefliger forms  $\Tr \left( A_t\exp(-(\delta A_t)^2\right)$ and $\Tr \left(  A'_t\exp (-(\delta A'_t)^2 \right)$ are closed.  In addition,
$$
\left[\Tr \left( A_t\exp(-(\delta A_t)^2) \right), \Tr \left( A'_t\exp(-(\delta A'_t)^2)\right) \right] \; \in \; H_c^*(M/F, M'/F'; \varphi)
$$
is independent of $t$.  So $ {{\ch}}(\Ind_a(D,D'))$ is well defined, and
$$
 {{\ch}}(\Ind_a(D,D')) \; = \; \Ind_t(D,D').
 $$
\end{theorem}
\begin{proof}
The Haefliger forms are closed by Proposition \ref{Cor3.7}, which also gives that $\frac{d}{dt}\Tr \left(  A_t\exp (-(\delta A_t)^2 \right) = d_H W_t$, and 
$\frac{d}{dt}\Tr \left(  A'_t\exp (-(\delta A'_t)^2 \right) = d_H W'_t$.  To finish the proof of $t$ independence, we need only show that $W_t$ and $W'_t$ can be chosen to be ${\varphi}$ compatible.

\ms
 
Recall that  $\pi_{\pm}$ is the grading operator 
$\pi_{\pm} \, = \, \diag(\Id_{\cS^+ \otimes E}, - \Id_{\cS^- \otimes E}) = \pi_+ - \pi_-$, and similarly for $\pi'_{\pm}$.  
When we identify the spin bundles and Dirac operators off compact subspaces, we also identify these gradings, so they are $\Phi$ compatible.  In particular, $\pi_-$ and $\pi'_-$ are $\Phi$ compatible.  Note that  $A_t + \pi_-$ and $A'_t + \pi'_-$ are idempotents.  Using this fact, in \cite{BH04}, Corrigendum, it is shown that 
$$
\frac{d}{dt}\left(\Tr \left( (A_t + \pi_-)\exp(-(\delta (A_t+ \pi_-))^2) \right), \Tr \left( (A'_t +  \pi'_-)\exp(-(\delta (A'_t +  \pi'_-)^2)\right) \right) \; = \; 
d_H(W_t,W'_t),
$$
where $(W_t,W'_t) \in \cA_c^*(M/F,M'/F'; \varphi)$, in particular they are ${\varphi}$ compatible.  This follows from the fact
 that the operators $\pa_{\nu}, \theta, \Theta,$ and $\delta$  all  preserve $\Phi$ compatability,  and that $W_t$ and $W'_t$ are constructed using those operators, $A_t$, $A_t'$, $\pi_-$, $\pi'_-$ {{(and the identities $\Id$ and $\Id'$)}}, and their derivatives with respect to $t$.  Since $A_t$, $A_t'$, $\pi_-$, $\pi'_-$, $\Id$ and $\Id'$ are $\Phi$ compatible, $W_t$ and $W'_t$ are $\varphi$ compatible.   As
 $$
d_H\Tr \left( (A_t + \pi_-)\exp(-(\delta (A_t + \pi_-))^2) \right) \; = \; d_H \Tr \left( (A'_t + \pi'_-)\exp(-(\delta (A'_t +  \pi'_-)^2)\right)  \; = \; 0,
$$
it follows that  
$$
\left[ \Tr \left( (A_t + \pi_-)\exp(-(\delta (A_t + \pi_-))^2) \right), \Tr \left( (A'_t + \pi'_-)\exp(-(\delta (A'_t +  \pi'_-)^2)\right)\right]\;  \; \in \; H_c^*(M/F, M'/F'; \varphi)
$$
is independent of $t$. 

\ms
 
Next, using Proposition 3.5 and Corollary 3.7 of \cite{BH08},  with the reasoning above, (that is: all the operators used in the proofs preserve $\Phi$ compatibility, so if the input is $\Phi$ compatible, the output is $\varphi$ compatible), shows that 
$$
\left[\Tr \left( A_t\exp(-(\delta A_t)^2) \right),\Tr \left( A'_t\exp(-(\delta A'_t)^2)\right) \right] \; = 
$$
$$
\left[ \Tr \left( (A_t + \pi_-)\exp(-(\delta (A_t+ \pi_-))^2) \right), \Tr \left( (A'_t + \pi'_-)\exp(-(\delta (A'_t + \pi'_-)^2)\right)\right] \;  \; \in \; H_c^*(M/F, M'/F'; \varphi).
$$

\ms

For the equality $\ch(\Ind_a(D,D')) = \Ind_t(D,D')$, standard techniques used in \cite{HL90, BH08}, coupled with Remark \ref{limseqs}, show that
$$
\lim_{t \to 0}\tr \left( A_t\exp(-(\delta A_t)^2) \right) \; = \; AS(D_F),
$$
uniformly pointwise on $M$, and we have the same for $A'_t$.  As $\Tr \left( A_t\exp(-(\delta A_t)^2) \right)$ involves  integrating over compact subsets,  we may interchange the limit with the integration.
\end{proof}

So we have Theorem  \ref{classicalindex}. 

\ms 

Note that so far, we have not used the assumptions in Theorem \ref{Secondhalf1} or  \ref{Secondhalf2}. We now move on to the proofs of Theorems \ref{Secondhalf1} and \ref{HaefTheorem}.
For the proof of Theorem \ref{Secondhalf1}, we need to show that 
$$
\lim_{t \to 1}\Tr \left( A_t\exp(-(\delta A_t)^2) \right) \; = \; \Tr \left( \pi_{\pm} P_0\exp(-(\delta (\pi_{\pm} P_0))^2) \right).
$$ 
Recall that $A_t$ is only defined for $0 < t <1$, that $P_0$ is the projection onto the kernel of $D^2$, and that $ A_t\exp(-(\delta A_t)^2)$  has propagation bounded by $c_A\sqrt{\beta(t)/\alp(t)}$ for some $c_A \in \R^+$.   
We recall below that  
$$
\lim_{t \to 1} k_{P_{[0,\ep)}A_tP_{[0,\ep)})}(\ovln{x},\ovln{x}) \; = \; k_{P_{[0,\ep)}  \pi_{\pm} P_0 P_{[0,\ep)}} (\ovln{x},\ovln{x}) \; =  \;  k_{\pi_{\pm} P_0}  (\ovln{x},\ovln{x}),
$$
uniformly pointwise which is sufficient for our purposes.  

\ms

Denote by $\varrho_{[\epsilon,\infty)}$ the characteristic function for the interval $[\epsilon,\infty)$.
\begin{lemma}\label{epPtest}  For $\ell$ a non-neqative integer, there exists a constant $C_{\ell}>0$ depending only on $\ell$,  such that
$$
\hspace{0.5cm}  ||z^{\ell}\varrho_{[\ep, \infty)}\chi^t(z)\varrho_{[\ep, \infty)}||_{\infty} \leq  C_{\ell}e^{-1/64\alp(t)} + \ep^{\ell}e^{- \ep^{2}\beta(t)/2}  \; \to \; 0, \text{  exponentially as } t \to 1.
$$
\end{lemma}
\begin{proof}
First note that, 
$$
||z^{\ell}\varrho_{[\ep, \infty)}\chi^t\varrho_{[\ep, \infty)}||_{\infty} \; \leq \; 
||z^{\ell}\varrho_{[\ep, \infty)}\left(\chi^t - e_{\sqrt{\beta (t)}}\right)\varrho_{[\ep, \infty)}||_{\infty} + 
||z^{\ell}\varrho_{[\ep, \infty)} e_{\sqrt{\beta (t)}}\varrho_{[\ep, \infty)} ||_{\infty} \; \leq \;
$$
$$        
||z^{\ell}\varrho_{[\ep, \infty)}\left(\chi^t - e_{\sqrt{\beta (t)}}\right)\varrho_{[\ep, \infty)}||_{\infty} + 
\ep^{\ell}e^{- \ep^{2}\beta(t)/2},
$$
since the maximum for the second term for $t$ close enough to $1$ will occur at $z = \ep$, as $\beta (t) \to \infty$ as $t \to 1$.  

\ms

Next, $||z^{\ell}\varrho_{[\ep, \infty)}\left(\chi^t - e_{\sqrt{\beta (t)}}\right)\varrho_{[\ep, \infty)}||_{\infty}$ is bounded by 
$||z^{\ell}\left(\chi^t - e_{\sqrt{\beta (t)}}\right)||_{\infty}$, which in turn is bounded by the $L^1$ norm of $FT(z^{\ell}\left(\chi^t - e_{\sqrt{\beta (t)}}\right))$.  Up to a constant depending only on $\ell$, 
$$
FT(z^{\ell}\left(\chi^t - e_{\sqrt{\beta (t)}}\right)) \; = \; \frac{\pa^{\ell}}{\pa z^{\ell}}FT\left(\chi^t - e_{\sqrt{\beta (t)}}\right)  \; = \; 
$$
$$
\frac{\pa^{\ell}}{\pa z^{\ell}} \left(\frac{1}{\sqrt{\beta(t)}} FT \left( \frac{1}{\sqrt{\alpha(t)}}\hat \psi_{1/\sqrt{\alpha (t)}} * e\right)_{1/\sqrt{\beta (t)}} - \frac{1}{\sqrt{\beta(t)}} e_{1/\sqrt{\beta(t)}} \right) \; = \;
$$
$$
\frac{\pa^{\ell}}{\pa z^{\ell}} \left(\frac{1}{\sqrt{\beta(t)}}\left(\psi_{\sqrt{\alpha (t)}}  e\right)_{1/\sqrt{\beta (t)}} \; - \; \frac{1}{\sqrt{\beta(t)}} e_{1/\sqrt{\beta(t)}} \right)\; = \; 
\frac{\pa^{\ell}}{\pa z^{\ell}} \left(\frac{e_{1/\sqrt{\beta (t)}}}{\sqrt{\beta(t)}} \Big(\psi_{\sqrt{\frac{\alpha (t)}{\beta(t)}}} - 1\Big)\right)   \; = \; 
$$
$$
\sum_{k = 0}^{\ell} {\left(^\ell_k\right)}\, 
\frac{\pa^k}{\pa z^k} \left(\frac{e_{1/\sqrt{\beta (t)}}}{\sqrt{\beta(t)}}\right) \frac{\pa^{\ell-k}}{\pa z^{\ell-k}}\Big(\psi_{\sqrt{\frac{\alpha (t)}{\beta(t)}}} - 1\Big).
$$
The function  $\psi_{\sqrt{\frac{\alpha (t)}{\beta(t)}}} - 1 = 0$ on  $\vert z\vert \leq \frac{1}{4\sqrt{\alpha (t)/\beta (t)}}$, and the norms of its derivatives are globally bounded by a constant depending only on $\ell$.  Thus the $L^1$ norm of {$FT\left(z^{\ell}\left(\chi^t - e_{\sqrt{\beta (t)}}\right)\right)$} is bounded by a constant, depending only on $\ell$, times 
$$
\sum_{k = 0}^{\ell}\int_{\vert z\vert \geq \frac{1}{4\sqrt{\alpha (t)/\beta (t)}} }
{\left\vert\frac{\pa^{k}}{\pa z^k} \left(e_{1/\sqrt{\beta (t)}}\right)\right\vert} \frac{dz}{\sqrt{\beta(t)}}    \; = \;
\sum_{k = 0}^{\ell}\int_{\vert z\vert \geq \frac{1}{4\sqrt{\alpha (t)/\beta (t)}} }
{\left\vert p_k(1/\sqrt{\beta} (t),z/\sqrt{\beta}) \right\vert}\left(e_{1/\sqrt{\beta (t)}}\right) \frac{dz}{\sqrt{\beta(t)}}    \; = \;
$$
$$
\sum_{k = 0}^{\ell}\int_{\vert z\vert \geq \frac{1}{4\sqrt{\alpha (t)}}} {\left\vert p_k(1/\sqrt{\beta} (t),z)\right\vert e^{-z^2/2}   dz \; \leq \; e^{-1/64\alpha (t)} \int_\R \sum_{k = 0}^{\ell} \left\vert p_k(1/\sqrt{\beta} (t),z)\right\vert e^{-z^2/4} dz}. 
$$
Here $p_k$ is a polynomial of degree $k$ in both variables, so the integral is bounded by a constant depending only on $\ell$.  Since $\alp(t) \to 0$ and $\beta(t) \to \infty$ as $ t \to 1$, we have the lemma.
\end{proof}
  
Denote by $Q_{\ep}$ the spectral projection for $D^{2}$ for the interval $[\epsilon,\infty)$, that is $Q_{\epsilon} = \varrho_{[\ep, \infty)} (D^2)$.  Since $I = P_{[0,\ep)} + Q_{\ep}$, the operators $Q_{\ep}$ and $\delta Q_{\ep}$ are transversely smooth and bounded, as the other two operators are because of our assumption of transverse smoothness.  The operators $P_{[0,\ep)}$, $Q_{\epsilon}$, and $A_t$ all commute as they are functions of $D$, so
$$
A_t\exp \left(-(\delta A_t)^2\right) = P_{[0,\ep)}A_tP_{[0,\ep)}\exp\left(-(\delta A_t)^2\right) \;+\; 
Q_{\ep}A_tQ_{\ep}\exp\left(-(\delta A_t)^2\right).
$$
Recall that $\delta A_t$ is nilpotent, in particular, 
$\dd \exp\left(-(\delta A_t)^2\right) = \sum_{k = 0}^{q/2} \frac{(-(\delta A_t)^2)^k}{k!}$, see \cite{BH08}.  Lemma \ref{epPtest} gives immediately that
$
||D^{2\ell}Q_{\ep}\chi^t(D)Q_{\ep}|| \to 0
$ 
exponentially as $t \to 1$.  The fact that every element of $A_t$ contains at least one $\chi^t(D)$, and that all the other terms are bounded, save $D^+$ (but $R_tD^+$  is covered by Lemma \ref{epPtest}), give that 
$
||D^{2\ell}Q_{\ep}A_tQ_{\ep}|| \to 0
$ 
exponentially as $t \to 1$.
Thus,
$
||D^{2\ell}Q_{\ep}A_tQ_{\ep}\exp\left(-(\delta A_t)^2\right) || \to 0
$ 
exponentially as $t \to 1$.  It follows from the proof of Theorem 2.3.9 and the statement of Theorem 2.3.13, both of \cite{HL90}, that the Schwartz kernel of $ Q_{\epsilon}A_tQ_{\epsilon}  \exp\left(-(\delta A_t)^2\right) \to 0$ pointwise uniformly exponentially as $ t \to 1$.  So,
$$
\lim_{t \to 1}\Tr \Bigl{(}  Q_{\epsilon}A_tQ_{\epsilon} \exp\left(-(\delta A_t)^2\right)\Bigr{)}\; = \;  0,
$$ 
in $\cA_c^*(M/F)$ 
and similarly for $Q'_{\epsilon}A'_tQ'_{\epsilon}$.   Thus we may ignore those terms.  Note carefully that this is true for {\em fixed} $\ep > 0$.

\ms

For the terms coming from $P_{[0,\ep)}A_tP_{[0,\ep)}$, note that for $t$ near $1$,
$2\varrho_{[0,\ep)}(z)$ dominates $\varrho_{[0,\ep)}\chi^t\varrho_{[0,\ep)}(z)$.  This follows from Remark \ref{limseqs}, since
$\lim_{t \to 1} \chi^t(z) = e^{-z^2/(2(1-t))}$ in the Schwartz topology, and for $t$ near $1$, $\sup_z e^{-z^2/(2(1-t))} = 1$.
Thus, for $\ell$ a fixed positive integer and for $t$ near $1$,
$$
2\ep^{\ell} ||\varrho_{[0,\ep})(z)||_{\infty}  \; = \;  2 ||z^{\ell}\varrho_{[0,\ep})(z)||_{\infty}  \; \geq \;  
||z^{\ell}\varrho_{[0,\ep)}\chi^t\varrho_{[0,\ep)}(z)||_{\infty}. 
$$
The fact that  $||\delta A_t||$ and  $||\exp \left(-(\delta A_t )^2\right)||$ are bounded and the argument above imply that a multiple of $\tr\left( k_{P_{[0,\ep)}}(\ovln{x},\ovln{x}) \right)$ dominates $||\tr \Bigl{(}  k_{P_{[0,\ep)}A_tP_{[0,\ep)}\exp\left(-(\delta A_t)^2\right)} (\ovln{x},\ovln{x})\Bigr{)} ||$.
Since we can ignore $Q_{\ep}A_tQ_{\ep}\exp\left(-(\delta A_t)^2\right)$, 
the Dominated Convergence Theorem implies 
$$
\lim_{t \to 1}\Tr \Bigl{(}  A_t\exp\left(-(\delta A_t)^2\right) \Bigr{)} \; = \; 
\lim_{t \to 1}\Tr \Bigl{(}   P_{[0,\ep)}A_tP_{[0,\ep)}
\exp\left(-(\delta A_t)^2\right) \Bigr{)}    \; = \; 
$$
$$
\lim_{t \to 1}\int_F \tr \Bigl{(}   P_{[0,\ep)}A_tP_{[0,\ep)}\exp \left(-(\delta A_t)^2\right) \Bigr{)}   \; = \; 
\int_F  \hspace{-0.0cm}\lim_{t \to 1}\tr \Bigl{(}   P_{[0,\ep)}A_tP_{[0,\ep)}\exp\left(-(\delta A_t)^2\right) \Bigr{)},
$$
and similarly for $A'_t$. 

\ms

The proof of Theorem 4.2 in \cite{BH08}, which requires that $F$ be Riemannian, shows that, under our conditions on the Novikov-Shubin invariants, in degree $2k$ for  $0 \leq 2k  \leq 2\ell$ we have, 
$$
\lim_{t \to 1}( P_{[0,\ep)}A_tP_{[0,\ep)}) \; = \; P_{[0,\ep)}  \pi_{\pm} P_0 P_{[0,\ep)} \; = \;  \pi_{\pm} P_0,
$$
uniformly pointwise, and similarly for $A'_t$.  So, in degree  $2k$ for $0 \leq 2k \leq 2\ell$,
$$
\int_F  \hspace{-0.0cm}\lim_{t \to 1}\tr \Bigl{(}   P_{[0,\ep)}A_tP_{[0,\ep)}\exp\left(-(\delta A_t)^2\right) \Bigr{)} \; = \; 
\int_F  \tr \Bigl{(}  \pi_{\pm}P_0\exp\left(-(\delta ( \pi_{\pm} P_0))^2\right)\Bigr{)} 
 \, = \,  \ch_{a}(P_{0}).
$$

So, we  have proven Theorems \ref{Secondhalf1} and \ref{HaefTheorem}.
 
\ms

It remains to prove Theorem \ref{Secondhalf2}, and we thus need to compute the limits as $t \to 0$ and $ t \to 1$ of
$$  
\lim_{s \to \infty} \left(\int_{T \ssm T_s}  \hspace{-0.5cm}\Tr \left( A_t\exp(-(\delta A_t)^2 \right) \wedge \ome_T \; - \;
 \int_{T' \ssm T_s'}  \hspace{-0.5cm}\Tr \left( A'_t\exp(-(\delta A'_t)^2\right) \wedge \ome'_{T'} \right).
$$
For $\lim_{t \to 0}$, we may assume that the two integrands agree on $M(0) = V_M$ and $M'(0)= V'_{M'}$ (actually on fixed penumbras).  Then we have, 
$$  
\lim_{t \to 0}\lim_{s \to \infty} \left(\int_{T \ssm T_s}  \hspace{-0.5cm}\Tr \left( A_t\exp(-(\delta A_t)^2 \right) \wedge \ome_T \; - \;
 \int_{T' \ssm T_s'}  \hspace{-0.5cm}\Tr \left( A'_t\exp(-(\delta A'_t)^2\right) \Bigr{)} \wedge \ome'_{T'} \right) 
 \; = \; 
$$
$$
\lim_{t \to 0}\left(\int_{M \ssm M(0)} \hspace{-0.7cm}\tr \left( A_t\exp(-(\delta A_t)^2 \right) \wedge \ome \; - \; 
 \int_{M' \ssm M'(0)} \hspace{-0.7cm}\tr \left( A'_t\exp(-(\delta A'_t)^2\right) \Bigr{)} \wedge \ome' \right) 
 \; = \; 
$$
$$
\int_{\cK} \hspace{-0.0cm} AS(D_F) \wedge \ome - 
\int_{\cK' } \hspace{-0.0cm}AS(D'_{F'})\wedge \ome'  \; = \; 
\left\langle(\left[  \int_F \AS(D_F),\int_{F'}AS(D_{F'}))  \right],(\ome_T,\ome'_{T'})\right\rangle. 
$$
As above, 
$
\lim_{t \to 0}\tr \left( A_t\exp(-(\delta A_t)^2 \right) \Bigr{)} \; = \; AS(D_F),
$
uniformly pointwise  on $M$, and we have the same for $A'_t$.  As we are integrating over compact subsets,  we may interchange the $\lim_{t \to 0}$ with the integrations.

\ms

For $\lim_{t \to 1}$, note that the operators have propagations bounded by $c_A\sqrt{\beta(t)/\alp(t)}$ for some $c_A \in \R^+$. As they are $\Phi$ compatible, we may assume that the two integrands agree on $T_{c_A(1-t)^{-1}}$ and $T'_{c_A(1-t)^{-1}}$.  Thus,
$$
 \lim_{t \to 1}\lim_{s \to \infty} \left(\int_{T \ssm T_s}  \hspace{-0.5cm}\Tr \left( A_t\exp(-(\delta A_t)^2 \right)\wedge \ome_T \; - \; 
 \int_{T' \ssm T_s'}  \hspace{-0.5cm}\Tr \left( A'_t\exp(-(\delta A'_t)^2\right) \Bigr{)} \wedge \ome'_{T'} \right)  \; = \; 
$$
$$
\lim_{t \to 1}\left(\int_{T \ssm T_{c_A(1-t)^{-1}}}  \hspace{-1.0cm}\Tr \left( A_t\exp(-(\delta A_t)^2 \right) \wedge \ome_T \; - \;
\int_{T' \ssm T'_{c_A(1-t)^{-1}}}  \hspace{-1.0cm}\Tr \left( A'_t\exp(-(\delta A'_t)^2\right) \Bigr{)} \wedge \ome'_{T'} \right). 
$$
Since the Schwartz kernel of $ Q_{\epsilon}A_tQ_{\epsilon}  \exp\left(-(\delta A_t)^2\right) \to 0$ pointwise exponentially as $ t \to 1$,  the  fact that $\ome_T$ is bounded, and the assumption that $M$ has sub-exponential growth, give that 
$$
\lim_{t \to 1} \left(\int_{T \ssm T_{c_A(1- t)^{-1}}}  \hspace{-1.0cm}\Tr \Bigl{(}  Q_{\epsilon}A_tQ_{\epsilon}  
\exp\left(-(\delta A_t)^2\right)\Bigr{)} \wedge \ome_T\right) \; = \;  0,
$$  
and similarly for $Q'_{\epsilon}A'_tQ'_{\epsilon}$.   Thus we may ignore those terms.

\ms 

Next, we have that for $t$ near $1$, a multiple of $\dd \int_M\tr\left( k_{P_{[0,\ep)}}\right)$ dominates
$$
||\int_M\tr \Bigl{(}  k_{P_{[0,\ep)}A_tP_{[0,\ep)}\exp\left(-(\delta A_t)^2\right)} \Bigr{)} ||.
$$
But this latter equals $\dd ||\int_M\tr \Bigl{(}  k_{{{P_{[0,\ep)}}}A_tP_{[0,\ep)}\exp\left(-(\delta A_t)^2\right)P_{[0,\ep)}} \Bigr{)}||$,
since $\dd \int_M \tr= \int_T  \int_F \tr = \int_T \Tr$, and $\Tr$ is a trace.  Thus, we need only show that a multiple of
$$
\tr\left(P_{[0,\ep)}(\ovln{x},\ovln{x}) \right) \; \text{dominates} \; ||\tr \Bigl{(} {{P_{[0,\ep)}}} A_tP_{[0,\ep)}\exp\left(-(\delta A_t)^2\right) P_{[0,\ep)}(\ovln{x},\ovln{x})\Bigr{)} ||.
$$
This is due to the fact that, for a smoothing operator $A$,
$
\tr(k_A((\ovln{x},\ovln{x})) = \sum_i \langle A ( \delta_{v_i}^{\ovln{x}}) , \delta_{v_i}^{\ovln{x}} \rangle.  
$
Here $v_i$ is an orthonormal basis of the fiber over the point $\ovln{x}$, and $ \delta_{v_i}^{\ovln{x}}$ is the Dirac delta section of the bundle supported at $\ovln{x}$.  Furthermore, everything is well defined on bounded geometry manifolds.
 See, for example, \cite{HL90} for details of such arguments.  As the operators we are concerned with are bounded leafwise smoothing,  we have,  
$$
||\langle {{P_{[0,\ep)}}} A_tP_{[0,\ep)}\exp\left(-(\delta A_t)^2\right)P_{[0,\ep)}(\delta_{v_i}^{\ovln{x}}),\delta_{v_i}^{\ovln{x}} \rangle|| { { \; = \; 
||\langle A_tP_{[0,\ep)}\exp\left(-(\delta A_t)^2\right)P_{[0,\ep)}(\delta_{v_i}^{\ovln{x}}),P_{[0,\ep)}(\delta_{v_i}^{\ovln{x}})\rangle||  }} \;\leq \;
$$
$$
||  A_tP_{[0,\ep)}\exp\left(-(\delta A_t)^2\right)|| \, {{||P_{[0,\ep)}(\delta_{v_i}^{\ovln{x}}) ||^2}}  \; = \;
||A_tP_{[0,\ep)}\exp\left(-(\delta A_t)^2\right)|| \,  \langle P_{[0,\ep)}(\delta_{v_i}^{\ovln{x}}), \delta_{v_i}^{\ovln{x}} \rangle. 
$$  
Summing over $i$, gives the result.

\ms

The fact that $||\ome||$ is bounded and the assumption that 
$\dd\int_M \tr\left( P_{[0,\ep)} \right)dx \; < \;  \infty$, imply that 
$$
\int_M \tr\left( P_{[0,\ep)} \right)||\ome|| \, dx 
\; < \;  \infty.
$$
Thus
$$
\int_M {{|| \tr \Bigl{(}  P_{[0,\ep)}A_tP_{[0,\ep)}\exp\left(-(\delta A_t)^2\right) (\ovln{x},\ovln{x})\Bigr{)} \wedge\ome||}} \, dx \; < \;  \infty,
$$
so {{the integral $\dd \int_M  \tr \Bigl{(}   P_{[0,\ep)}A_tP_{[0,\ep)}\exp\left(-(\delta A_t)^2\right) \Bigr{)} \wedge \ome$ converges. Notice that
$$
\int_M  \tr \Bigl{(}   P_{[0,\ep)}A_tP_{[0,\ep)}\exp\left(-(\delta A_t)^2\right) \Bigr{)} \wedge \ome \; = \; 
\int_T  \Tr \Bigl{(}   P_{[0,\ep)}A_tP_{[0,\ep)}\exp\left(-(\delta A_t)^2\right) \Bigr{)} \wedge \ome_T.
$$
This fact, the fact that we can ignore $Q_{\ep}A_tQ_{\ep}\exp\left(-(\delta A_t)^2\right)$, 
and the Dominated Convergence Theorem imply, 
$$
\lim_{t \to 1}\int_{T \ssm T_{c_A(1- t)^{-1}}}  \hspace{-1.5cm}\Tr \Bigl{(}  A_t\exp\left(-(\delta A_t)^2\right) \Bigr{)} \wedge \ome_T \; = \; 
\lim_{t \to 1}\int_{T \ssm T_{c_A(1-t)^{-1}}}  \hspace{-1.5cm}\Tr \Bigl{(}   P_{[0,\ep)}A_tP_{[0,\ep)}
\exp\left(-(\delta A_t)^2\right) \Bigr{)} \wedge \ome_T   \; = \; 
$$
$$
\lim_{t \to 1}\int_{T} \Tr \Bigl{(}   P_{[0,\ep)}A_tP_{[0,\ep)}\exp \left(-(\delta A_t)^2\right) \Bigr{)} \wedge \ome_T   \; = \; 
\int_T  \hspace{-0.0cm}\lim_{t \to 1}\Tr \Bigl{(}   P_{[0,\ep)}A_tP_{[0,\ep)}\exp\left(-(\delta A_t)^2\right) \Bigr{)} \wedge \ome_T,
$$
and similarly for $A'_t$.  The proof of Theorem 4.2 in \cite{BH08} shows that, as above, under our conditions on the Novikov-Shubin invariants in Theorem \ref{Secondhalf2}, 
$$
\lim_{t \to 1}\Tr \Bigl{(}  A_t\exp\left(-(\delta A_t)^2\right)\Bigr{)} \wedge \ome_T \; = \; 
\Tr \Bigl{(} \pi_{\pm} P_0\exp\left(-(\delta( \pi_{\pm}  P_0))^2\right)\Bigr{)} \wedge \ome_T, 
$$
and similarly for $A'_t$.  So,
$$
\lim_{t \to 1}\left(\int_{T \ssm T_{c_A(1-t)^{-1}}}  \hspace{-1.5cm}\Tr \Bigl{(}  A_t\exp\left(-(\delta A_t)^2\right)\Bigr{)} \wedge \ome_T - \int_{T' \ssm T'_{c_A(1-t)^{-1}}}  \hspace{-1.5cm}\Tr \Bigl{(}  A'_t\exp \left(-(\delta A'_t)^2\right) \Bigr{)} \wedge \ome'_{T'} \right) \; = \; 
$$
$$
\int_T  \Tr \Bigl{(}  \pi_{\pm}  P_0\exp\left(-(\delta ( \pi_{\pm}  P_0))^2\right)\Bigr{)} \wedge \ome_T   \; - \; 
\int_{T'}  \Tr \Bigl{(}  \pi_{\pm}' P'_0\exp\left(-(\delta( \pi_{\pm}' P'_0))^2\right) \Bigr{)} \wedge \ome'_{T'}  
 \, = \,
 $$
 $$   
{{\langle(\ch_{a}(P_{0}),\ch_{a}(P'_{0})),{{(\ome_T,\ome'_{T'})}}\rangle.}}
$$
That is,   
$
 \langle {{{\ch}} \Ind_a(D,D'), {{ [\ome_T,\ome'_{T'}]}}}\rangle   =  \langle(\ch_{a}(P_{0}),\ch_{a}(P'_{0})),{{(\ome_T,\ome'_{T'})}}\rangle.
$
So we have proven Theorem \ref{Secondhalf2}. 

\begin{remark} \label{growth} 
Note that if $M$, so also $M'$, grows exponentially, there are constants $c_0, c_M \in \R^+$, so that $\vol(M_t) \leq c_0e^{c_M t}$.  This follows from the Bishop-Gromov inequality.  Thus, if we used Lemma \ref{epPtest} as above and integrated over $M \ssm M_{c_A\sqrt{\beta(t)/\alp(t)}}$, we would get an estimate of the form,
$$
(C_{\ell}e^{-1/64\alp(t)} + \ep^{\ell}e^{- \ep^{2}\beta(t)/2})c_0e^{c_Mc_A\sqrt{\beta(t)/\alp(t)}}.
$$
For the proof to work, we need this to $\to 0$ as $t \to 1$.  Now  $\sqrt{\beta(t)/\alp(t)} = \beta(t) = 1/\alp(t)$, as $t \to 1$.
Thus the two terms must $\to 0$ individually.  This only happens if  $c_Mc_A <  \min(\ep^2/2, \frac{1}{64})$.  That is, the exponential growth is not too robust.
\end{remark}

\section{Invertible near infinity operators}\label{TIMAFTS}

In this section, we assume that $(M,F)$ is as in the first two paragraphs of Section \ref{setup}.  

\subsection{Invertibility near infinity}

Our new assumption here is that the zero-th order contribution $\maR^E_F$ in the Bochner formula defined below is strictly positive on $M$ near infinity.  As $\maR^E_F$ is locally defined, this implies that the same for $\maR^{E'}_{F'}$.

\ms

For the leafwise Dirac operator $D_F = (D_L)$,  the canonical operator $\maR^E_F$ on sections of $E_M|_L$ is given by
$$
\maR^E_F(\varphi) \,\, = \,\, \frac{1}{2} \sum_{i,j=1}^p \rho(X_i) \rho(X_j)  R^E_{X_i,X_j} (\varphi),
$$ 
where $R^E$is the curvature operator of the Hermitian connection $\nabla^{F,E}$ on $E_M|_L$, 
$X_1,\ldots, X_p$ is a local oriented orthonormal basis of $TL$, and  $\rho(X_i)$ is the Clifford action of $X_i$. 
Note that $\maR^E_F$ is well defined, smooth, and that it is globally bounded because of the assumption  of bounded geometry.   
The operators $D_L$  and $\maR^E_F$ are related by the general leafwise Bochner Identity, \cite{LM89}
\begin{Equation}\label{basic}
$
\hspace{4.9cm} D_L^2 \,\, = \,\, (\nabla^{F,E})^* \nabla^{F,E} \,\, + \,\,  \maR^E_F.
$ 
\end{Equation}
As we work on $\cG$ rather than $M$,  $D = r^*(D_F)$ also satisfies Equation \ref{basic}, which, being local, is  the same, namely, $D^2 \,\, = \,\, \nabla^* \nabla \,\, + \,\,   r^*(\maR^E_F)$.  Note that in general, if $\maR^E_F$ is strictly positive near infinity, $r^*(\maR^E_F)$ is not, due to the fact that {{$r$ is not a proper map in general}}.  However, $r^*(\maR^E_F)$ is $\maG$-invariant strictly positive near infinity off some $\maG$-compact subspace, in particular when restricted to $M \subset \cG$, since it coincides with $\maR^E_F$ there.  

\ms

 We have the following result from of \cite{BH21}.   Note that it does not need $P_{[0, \ep]}$ to be transversely smooth.  It does need it to be transversely measurable, which it is by Lemma 4.10 of \cite{BH21}.  
\begin{theorem}\label{measPRedux} {(Theorem 5.2 of \cite{BH21})}  Assume that $F$ admits a holonomy invariant transverse measure $\Lam$.  Suppose $\maR^E_F$  is strictly positive near infinity.  In particular, we may assume that $\kappa_0 = \sup \{ \kappa \in \R \, | \, \maR^E_F -\kappa \Id \geq 0 \; \text{on} \; M\smallsetminus \cK_M\}$ is positive.    Then, for $0 \leq \ep  < \kappa_0$, 
$$
 \int_M   \tr(P_{[0,\ep]}(\overline{x},\overline{x}))  \, dx_F d\Lam  \;  \leq \; 
\frac{(\kappa_0 - \kappa_1)}{(\kappa_0 - \ep)}   \int_{\cK_M}   \tr(P_{[0,\ep]}(\overline{x},\overline{x})) \, dx_F d\Lam \,\, < \,\, \infty,   
$$
where $\kappa_1  = \sup \{ \kappa \in \R \, | \, \maR^E_F -\kappa \Id \geq 0 \; \text{on} \; M\}$.    
\end{theorem} 

Note that if $F$ is Riemannian, it does admit holonomy invariant transverse measures, and we can insure that $dx$ is of the form $dx_F d\Lam$. 

\begin{proof}
The proof of Theorem 5.2 in \cite{BH21} works equally well here, mutatis mutandis.   The changes in notation needed are 
$$
D^E_L \to D, \; k_{P_{[0,\ep]}}(x,x) \to k_{P_{[0,\ep]}}(\ovln{x},\ovln{x}), \; L \to \wtit{L}_x, \; \sigma_L \to \sigma_x,\;  \int_L \to \int_{\wtit{L}_x},
$$
and so on.  
\end{proof}

Proposition 5.5 in \cite{BH21} still holds here, namely the following.
\begin{prop}\label{measPR2} Suppose the curvature operator $\maR_F^E$ is strictly positive on $M$, that is $\kappa_1 >0$, so $\maR_F^E\geq \kappa_1 \Id$ on $M$.   Then for $0 \leq \ep < \kappa_1$, 
$P_{[0,\ep]} = 0.$
\end{prop}

The relationship with the index bundle is not insured in general, \cite{BHW14}, and one needs to impose additional  spectral assumptions.
We have, as in \cite{BH21}, the following immediate corollaries of Theorems  \ref{Secondhalf1}, \ref{Secondhalf2} and \ref{measPRedux} which relate the pairings there to pairings with the index bundles. 

\begin{theorem}\label{IMCOR1}
Suppose that  $(M, F, \cK_M)$ and $(M', F',\cK_{M'})$ are bounded geometry foliations which are identified outside the compact subspaces $\cK_M$ and $\cK_{M'}$ as before and let $(\omega, \omega')$ be a $\varphi$-compatible pair of closed holonomy invariant forms of degree $\ell\leq q$.  Assume the following:
\begin{itemize}
\item $M$, and so also $M'$, has sub-exponential growth, and $F$ and $F'$ are Riemannian;
\item the leafwise operators $P_0$, $P_0'$, $P_{(0,\ep)}$  and $P_{(0,\ep)}'$ (for $\epsilon$ sufficiently small) 
are transversely smooth;
\item $NS(D)$ and $NS(D')$ are greater than {{$\ell$}};
\item $\maR^E_F$, so also $\maR^{E'}_{F'}$,  is strictly positive near infinity in $M$ and $M'$ respectively. 
\end{itemize}
Then 
$$
\langle\left[\int_F AS(D_F),\int_{F'} AS(D_{F'})\right],[\ome_T,\ome'_{T'}]\rangle \, = \,   
\langle(\ch(P_{0}),\ch(P'_{0}),(\ome_T,\ome'_{T'})\rangle.
$$
\end{theorem}

Recall that  $(\ch_{a}(P_{0}),\ch_{a}(P'_{0})$ is not an element of $\maA_c^*(M/F, M'/F'; \varphi)$ in general.
Since $AS(D_F)$ and $AS(D_F)$ are ${\varphi}$ compatible, $AS(D_F) \wedge \ome$ and  $AS(D_{F'}') \wedge \ome$ are ${\varphi}$ compatible, say off the compact subsets $\what{\cK}$ and $\what{\cK}'$, and then we have 
$$
\langle\left[\int_F AS(D_F),\int_{F'} AS(D_{F'})\right], [\ome_T,\ome'_{T'}]\rangle \, = \,  
\int_{\what{\cK}} AS(D_F) \wedge \ome \, -  \, \int_{\what{\cK}'} AS(D_{F'}) \wedge \ome'.
$$
  
\ms

For a single foliated manifold we have the following, compare with \cite{GL83}.
 
 \begin{theorem}\label{IMCOR2} 
Suppose that $E$ and $E'$ are two Clifford bundles over the foliated manifold $(M,F)$, which are isomorphic off the compact subset $\cK_M$, with associated Dirac operators $D$ and $D'$.  Let $\omega$ be a bounded closed holonomy invariant transverse form (or current) of degree $\ell\leq q$.  Suppose that
\begin{itemize}
\item $M$ has sub-exponential growth, and $F$ is Riemannian;
\item the leafwise operators $P_0$, $P_0'$, $P_{(0,\ep)}$  and $P_{(0,\ep)}'$ (for $\epsilon$ sufficiently small)  are transversely smooth;
\item $\min (NS(D), NS(D'))$ is greater than $\ell$;
\item $\maR^E_F$, and hence also $\maR^{E'}_F$,  is strictly positive near infinity.   
\end{itemize}
 Then, since $\ch(E)  = \ch(E')$ off $\cK_M$,
$$
\langle\ch(\Ind_a(D,D')), [\ome_T,\ome_T]\rangle \, = \,\int_{\cK_M} (\AS (D_F)(\ch(E) - \ch(E')) \wedge \ome 
\, = \,  \langle(\ch(P_{0}),\ch(P'_{0})),(\ome_T,\ome_T)\rangle.
$$
\end{theorem}

\begin{remark}
Note that if $E_1$ is a leafwise almost flat bundle (actually a $K$-theory class) on $M$, then we may twist the operators $D$ and $D'$ by $E_1$ to get the operators $D_{E_1}$ and $D'_{E_1}$.  Uniform positivity near infinity is preserved when this is done, so we have the extension of Theorem \ref{IMCOR2}  to $D_{E_1}$ and $D'_{E_1}$.  Theorem \ref{IMCOR1} also extends in this way if we have  leafwise almost flat bundles  $E_1 \to M$ and $E_1' \to M'$ which are isomorphic near infinity.
\end{remark}

\subsection{Reflective foliations}

We now relate our definition of the relative index to the cut-and-paste definition considered in Section 4 of \cite{GL83}.  For this paper to be self-contained, we paraphrase from \cite{BH21}.  For simplicity, we assume that $\ome$ and $\ome'$ are ${\varphi}$ compatible off $\cK_M$ and $\cK'_{M'}$.

\ms

We say that $(M, F, \cK_M)$ as above is reflective if there exists a compact submanifold $H\subset M$  such that
$$
\cK_M\subset H\; \text{ and }\; \pa H\text{ is transverse to }F.
$$
So $F'$ is also reflective with corresponding $H'$.  Then there is $\delta >0$, and a neighborhood  of {{$\pa H$ which is diffeomorphic to $\pa H \times [-\delta, \delta]$, and so that $F$ restricted to $\pa H \times [-\delta, \delta]$ has leaves of the form $(L \cap \pa H) \times[-\delta, \delta]$.  We may assume that the foliation preserving diffeomorphism $\varphi$ extends to $\pa H\times  [-\delta,\delta]$, and  that  $\varphi(\pa H \times [-\delta, \delta])$ is diffeomorphic to $\pa H' \times [-\delta, \delta]$, and that it has the same properties as $\pa H \times [-\delta, \delta]$.  Then we may form the compact foliated manifold 
$$
\what{M}  \,\,= \,\, H \cup_{\what{\varphi}} H',
$$
where $\what{\varphi}:\pa H \times [-\delta , \delta] \to \pa H' \times [-\delta, \delta]$ is given   by   $\what{\varphi} (x,s) = (\varphi (x),-s)$}}.  We change the orientation of $F'$ to the opposite of what it was originally.  The resulting foliation $F \cup_{\what{\varphi}} F'$ is denoted $\what{F}$.  Denote by {{$\pi: \pa H \times [-\delta, \delta] \to \pa H$ the projection and   note that  $E \, |_{ \pa H \times [-\delta, \delta]} \simeq \pi^*(E \, |_{\pa H})$, and  $TF \, |_{ \pa H \times [-\delta, \delta]} \simeq \pi^*(TF \, |_{\pa H})$.  (Note that    $\dim(TF \, |_{\pa H}) = \dim (TF)$, not $\dim (TF) -1 = \dim F \, |_{\pa H}$.) }}  We may assume that  $\nabla$ and $D_F$ are preserved under the maps $(x,s) \to (x,-s)$ and $E_{(x,s)} \to E_{(x,-s)}$. This implies that $D_F$ and $D'_{F'}$ are identified under the  glueing map.  In addition,  $\ome$ and $\ome'$ fit together, giving $\what{\ome}$. This construction is the exact translation of the Gromov-Lawson construction to foliations.

\ms

Finally, denote the leafwise operator on $\what{F}$ by $\what{D}_{\what{F}}$ (and its associated projections by $\what{P}_0$ and  $\what{P}_{(0, \ep)}$).  Then we have the following extension of Alain Connes' celebrated index theorem, see \cite{C79}, which is very useful in Section \ref{apps}.

 \begin{theorem}\label{ExRelInd2}  With the above notations, suppose that  $F$ (and so also $F'$)  is reflective, but not necessarily Riemannian.  Suppose further that $(C,C')$ is a compatible near infinity pair of closed holonomy invariant currents, with associated current $\what{C}$.
Then   
$$
\langle \ch ( \Ind_a(D, D'), (C,C') \rangle  \: =  \: \langle \ch ( \Ind_a(\what{D}_{\what{F}})), \what{C} \rangle. 
$$
\end{theorem}

\begin{proof}
We prove the case where $(C,C') = [\ome_T,\ome'_{T'}]$, since it is notationally simpler.  The general case is left to the reader.
Theorems \ref{classicalindex} and \ref{Hformsclosed} give
$$
\langle \ch (\Ind_a(D, D')),[\ome_T,\ome'_{T'}]\rangle  \; =  \; \int_{\cK_M} \AS (D_F) \wedge \ome  -   \int_{\cK'_{M'}} \AS (D'_{F'})  \wedge \ome',
$$
since the differential forms $\AS (D_F) \wedge \ome$ and $\AS (D'_{F'})  \wedge \ome'$ are $\varphi$ compatible off $\cK_M$ and $\cK'_{M'}$.  Next, 
$$
\int_{\cK_M} \AS (D_F) \wedge \ome  -   \int_{\cK'_{M'}} \AS (D'_{F'})  \wedge \ome'  \; =  \;  \int_H \AS (D_F) \wedge \ome  -   \int_{H'} \AS (D'_{F'})  \wedge \ome' \; =  \; 
$$
$$
\int_{\what{M}}  \AS(\what{D}_{\what{F}}) \wedge \what{\ome} \; =  \;  \langle \int_{\what{F}}  \AS(\what{D}_{\what{F}}), \what{\ome} _{\what{T}}\rangle  \: =  \: \langle \ch ( \Ind_a(\what{D}_{\what{F}})), \what{\ome}_{\what{T}} \rangle.
$$
The last equality is from Theorem 6.2 of \cite{BH04} applied to the closed foliated manifold $(\what M, \what F)$.  The others are obvious.
 \end{proof}

Note that, since the integrands $\AS (D_F) \wedge \ome$ and $\AS (D'_{F'})  \wedge \ome'$ are ${{\varphi}}$ compatible off $\cK_M$ and $\cK'_{M'}$, this result is actually independent of the choice of the transverse compact hypersurface {$\pa H$ and for simplicity we may assume that $H = \cK_M$. 

\begin{theorem}\label{ExRelInd3}
Suppose  that $(M,F,D)$, $(M',F',D')$ are as in Theorem \ref{ExRelInd2}. Suppose furthermore that $\what F$ is Riemannian, that $\what{P}_0$ and $\what{P}_{(0, \ep)}$ are transversely smooth, and the Novikov-Shubin invariants of $\what{D}_{\what{F}}$ are greater than $\ell/2$, for some $0\leq \ell \leq q$.  Then for any $\ell$ homogeneous $\varphi$-compatible pair $(\ome, \ome')$ as before, 
$$ 
\langle \ch( \Ind_a(D,D')), [\ome_T,\ome'_{T'}]\rangle  =  \langle(\ch(\what{P}_{0}),\what{\ome}_{\what{T}}\rangle.
$$
Moreover, if we impose on $(M, F, D)$ and $(M', F', D')$ the assumptions of Theorem \ref{Secondhalf2}, then we have
$$
\langle (\ch(P_0), \ch(P_0)),  (\ome_T, \ome'_{T'})  \rangle \; = \; \langle(\ch (\what{P}_{0}),\what{\ome}_{\what{T}}\rangle.
$$
\end{theorem}

This is a consequence of Theorem \ref{ExRelInd2} using Theorem 4.1 of \cite{BH08} to deduce the second  equality, with $\langle(\ch(\what{P}_{0}),\what{\ome}_{\what{T}}\rangle$ being well defined under our assumptions.

\begin{remark}\ This result raises some interesting questions.  
\begin{enumerate}
\item  Suppose that $\maR_F$, so also $\maR_{F'}$, is strictly positive near infinity, then $\langle \ch(P_{0}),\ome_T \rangle$ and $\langle \ch(P'_{0})),\ome'_{T'} \rangle$ exist.  Under what more general conditions than those in Theorems \ref{Secondhalf2} and  \ref{ExRelInd3} does 
$$
\langle \ch(P_{0}),\ome_T \rangle - \langle \ch(P'_{0}),\ome'_{T'}\rangle =
\langle\ch(\what{P}_{0}),\what{\ome}_{\what{T}}\rangle?
$$
\item In general,  suppose that $\langle \ch(P_{0}),\ome_T \rangle - \langle \ch(P'_{0}),\ome'_{T'}\rangle -
\langle\ch(\what{P}_{0}),\what{\ome}_{\what{T}}\rangle \neq 0.$
What can be said about the geometry or topology of $(M,F,D)$, $(M',F',D')$, and $(\what{M}, \what{F},\what{D})$?
\item
 How are the Novikov-Shubin invariants of $D$ and $D'$ related to those of $\what{D}$?
\end{enumerate}
\end{remark}

The previous construction  extends to the following more general situation to yield the so called {{higher}} $\Phi$ relative index theorem, see again \cite{GL83}.  In particular, we assume that $(M,F)$ and $(M',F')$ satisfy the hypotheses of Theorem \ref{Secondhalf2}, with the following changes.  In particular,  $M \ssm \cK = V_+ \cup V_{\Phi}$  and  $M' \ssm \cK'  = V_+' \cup V_{\Phi}'$, where the unions are disjoint.  For this case,  $\Phi = (\phi, \varphi)$ is a  bundle morphism from $E\to V_{\Phi}$ to $E'\to V_{\Phi}'$ as in Section \ref{setup}, our good covers $\cU$ and $\cU'$ are compatible on $V_{\Phi}$ and $V_{\Phi}'$, and $\ome$ and $\ome'$ are $\varphi$ compatible on $V_{\Phi}$ and $V'_{\Phi}$.  We assume that $F$ is transverse to $\pa V_{\Phi}$, so $F'$ is transverse to  $\pa V_{\Phi}'$.  Finally, we assume  that  $\cR_F^E$ and $ \cR_{F'}^{E'}$ are strictly positive off $\cK$ and $\cK'$, so we do not need the assumptions on the integrals being finite.  

\ms

Next, consider as above the manifold
$\what{M} = (M \ssm V_{\Phi}) \cup_{\what{\varphi}} (M' \ssm V'_{\Phi})$, with the foliation
$$
\what{F} \; = \;  (F |_{M \ssm V_{\Phi}}) \cup_{\what{\varphi}} (F' |_{M' \ssm V'_{\Phi}}), 
$$ 
where the orientation on $\what{F}  |_{M \ssm V_{\Phi}} $ is the one on $F$, and that on $\what{F} |_{M' \ssm V'_{\Phi}} $ is the opposite of the  one on $F'$.  We also have the bundle $\what{E} \to \what{M}$ induced by $E$ and $E'$, the leafwise operator $\what{D}_{\what{F}}$ induced by $D_F$  and $D'_{F'}$,  and the differential form $\what{\ome}$ induced by $\ome$ and $\ome'$.

\ms

Because of the positivity off compact subsets, all three operators have finite indices, and we have the following.
\begin{theorem}\label{RelIndThm2}  [The higher foliated $\Phi$-index formula]
$$
\langle\ch( \Ind_a(\what{D})), \what{\ome}_{\what{T}} \rangle \, = \,  
\langle \ch(\Ind_a(D)), \ome_T \rangle \,\, - \,\, \langle \ch( \Ind_a(D')), \ome'_{T'} \rangle.
$$
\end{theorem} 

The proof  follows from our results here, by easily adapting the proof of Theorem 4.35 of  \cite{GL83}.

\section{Applications}\label{apps}

\subsection{Leafwise PSC and the higher Gromov-Lawson invariant}

We further extend the Gromov-Lawson construction in \cite{GL83}, Section 3, see also \cite{LM89}, IV.7,  to get an invariant for the  space of  PSC metrics on a foliation $F$ whose tangent bundle admits a spin structure.    We calculate this invariant for a large collection of spin foliations whose Haefliger $\what{A}$ genus is zero, so the results of \cite{BH21} do not apply.  Using the higher index results here, we show that the space of  PSC metrics on each of these foliations has infinitely many path connected components, thus verifying our claims that higher order index theorems allow for the extension of results for manifolds with non-zero $\what{A}$ genus to arbitrary manifolds, and that the higher order terms of the $\what{A}$ genus also carry geometric information.

\ms
For simplicity, we assume that $M$ is compact.  Denote by $\cM$ the space of all smooth metrics on  $F$ with the $C^\infty$ topology, and  by $\cM_{sc}^+  \subset \cM$ the subspace of  metrics with PSC along the leaves.  

\ms
Scalar curvature and the so called Atiyah-Singer operator are intimately related.  Recall that  $\maS_F$ is the bundle of  spinors along the leaves of $F$, with the leafwise spin connection $\nabla^F$.  The leafwise Atiyah-Singer operator is the leafwise spin Dirac operator $D_F^{\cS} = (D_L^\cS)$, which acts on $\cS_F$, as usual, by  
$$
D_L^{\cS} \; = \; \sum^p_{i=1}\rho(X_i)\nabla^F_{X_i},
$$
where $X_1,\ldots, X_p$ is a local oriented orthonormal basis of $TL$, and 
$\rho(X_i)$ is the Clifford action of $X_i$ on the bundle ${\cS_F} | L$.
Denote by $\kappa$ the leafwise scalar curvature of $F$, that is 
$$
\kappa = -\sum_{i,j =1}^p \langle R_{X_i,X_j}(X_i),X_j \rangle,
$$
where $R$ is the curvature operator associated to the metric on the leaves of $F$.  
In this case the Bochner Identity, Equation \ref{basic}, is quite simple, see \cite{LM89}, namely
\begin{Equation}\label{basicLisch}
$
\hspace{4.9cm} (D_L^{\cS})^2 \,\, = \,\, (\nabla^F)^* \nabla^F \,\, + \,\,  \frac{1}{4}\kappa.
$ 
\end{Equation}
Consider the foliation $F_{\R}$ on $M_{\R} = M \times \R$ with leaves $L_{\R} = L \times \R$ and with the  leafwise volume form $dx_F \times dt$.  If $\maU$ is a good cover of $M$,  $\cU_{\R} = \{(U_i^n,T_i^n) = (U_i \times (3n-2,3n+2), T_i) | \, (U_i,T_i) \in \maU, n \in \Z \}$ is a good cover of $M_{\R}$.  
Denote by $\pi:M_{\R} \to M$ the projection.  Suppose that $g_0, g_1 \in \cM_{sc}^+$,   and $(g_t)_{t \in [0,1]}$ is a smooth family in  $\cM$ from   $g_0$ to $g_1$.  On $F_{\R}$, set  $G = g_0 + dt^2$ for $t \leq 0$, $G = g_1 + dt^2 $ for $t \geq 1$, and $G = g_t + dt^2$ for $0 < t < 1$.  

\ms

The leafwise spin Dirac operator $D_F^{\cS}$ extends to the leafwise spin Dirac operator $D_{\R}$ on $F_{\R}$.   Following Gromov-Lawson, \cite{GL83}, Equation (3.13), we set
\begin{Equation}\label{GLinvar}
\hspace{4.5cm}$i(g_0,g_1) \; = \;    \ch(\Ind_a(D_{\R})) \;   \in  \; H^*_c(M_{\R}/F_{\R}).$
\end{Equation}

\begin{theorem} \label{Phithm} $i(g_0,g_1)$  depends only on $g_0$ and $g_1$.
If $i(g_0,g_1) \,\, \neq \,\, 0$, then $g_0$ and $g_1$ are not in the same path connected component of $\cM_{sc}^+$.  
\end{theorem}

\begin{proof}
 Suppose that $g_t$ and $\what{g}_t$ are two smooth families of metrics in $ \cM$ from $g_0$ to $g_1$, with associated metrics $G$ and $\what{G}$  and associated operators $D_{\R}$ and $\what{D}_{\R}$.  A byproduct of  Theorem \ref{classicalindex} is that $\dd i(g_0,g_1) \; = \; \left[ \int_{F_{\R}} \what{A}(TF_{\R}, {G})\right]$, where $\what{A}(TF_{\R}, {G})$ is the Atiyah-Singer characteristic differential form, the so-called  $A$-hat form of $F$, on $M_{\R}$ associated to the metric $G$, and similarly  for  $\what{G}$.  Thus we have
$$
i(g_0,g_1)(G)  -  i(g_0,g_1)(\what{G}) = 
\left[\int_{F_{\R}}  \hspace{-0.2cm}  (\what{A}(TF_{\R},{G}) \, -  \, \what{A}(TF_{\R}, {\what{G}})) \right]  \; = \; 
\left[\int_{F_{\R}^0}   \hspace{-0.2cm}  (\what{A}(TF_{\R},{G})  -  \what{A}(TF_{\R},{\what{G}}) )  \right], 
$$
where  $F_{\R}^0$ is the foliation on $\bigcup_i U_i^0$.   
 The forms $\what{A}(TF_{\R}, {G})$ and  $\what{A}(TF_{\R},{\what{G}})$ are locally computable in terms of their associated curvatures.  Thus, off the compact subset $M \times [0,1]$, they agree, which justifies the second equality.  By abuse of notation, we may write 
$$
\int_{F_{\R}^0}   \hspace{-0.1cm}  (\what{A}(TF_{\R},{G})  -  \what{A}(TF_{\R}, {\what{G}}))\; = \; 
\int_{F \times [0,1]}   \hspace{-0.6cm}  (\what{A}(TF_{\R},{G})  -  \what{A}(TF_{\R}, {\what{G}})). 
$$
Since the cohomology classes of the two forms are the same, $\what{A}(TF_{\R},{G})  -  \what{A}(TF_{\R},{\what{G}})$ is an exact form $d_{M \times\R}\Psi$, which is locally computable in terms of the curvatures and connections.   In particular, $\Psi = 0$ on the closure of open sets where the connections agree.   So off $M \times (0,1)$, $\Psi$  is zero, since the connections agree there.   Thus 
$$
\int_{F_{\R}^0}   \hspace{-0.1cm}  (\what{A}(TF_{\R},{G})  -  \what{A}(TF_{\R}, {\what{G}}))\; = \; 
\int_{F \times [0,1]}   \hspace{-0.6cm}  d_{M \times\R}\Psi \; = \; 
d_{H}\int_{F \times [0,1]}   \hspace{-0.6cm}  \Psi, 
$$
and $i(g_0,g_1)(G)  \,\, - \,\, i(g_0,g_1)(\what{G}) = 0$  in  $H^*_c(M_{\R}/F_{\R})$.

\ms

For the second part, assume that $g_0$ and $g_1$ are in the same path connected component of $\cM_{sc}^+$, and that $g_t$, is a  smooth family of metrics in $\cM_{sc}^+$ from $g_0$ to $g_1$.   Then $G$ restricted to each leaf of $F_{\R}$ has PSC, and since the family of metrics is smooth, it is strictly positive.  Then, Proposition \ref{measPR2} gives that  $P_{[0,\ep]}^{\R} = 0$, for some positive $\ep$, so the Novikov -Shubin invariants are infinite and Remark \ref{remarks4.4} (1) gives that $i(g_0,g_1) = 0$.   
\end{proof}

\begin{remark}
Theorem \ref{Phithm}  remains true if we consider concordance classes of PSC metrics, which a priori is stronger.  Recall that leafwise metrics are concordant if there is a metric $G$ on $TF_{\R}$ so that it agrees with $g_0 \times dt^2$ near $-\infty$ and with $g_1 \times dt^2$ near $+\infty$.  The conclusion is that if  $i(g_0,g_1) \,\, \neq \,\, 0$, then $g_0$ and $g_1$ are not in the same concordance class of metrics in $\cM^+_{sc}$.  The proof being essentially the same. 
\end{remark}

\begin{remark}
We could also extend this theory to concordance classes of leafwise flat connections $\nabla$ on an auxiliary bundle $E$.  The invariant would become $i((g_0, \nabla_0),(g_1, \nabla_1))$. See \cite{Be20}.  The theorem would then be that if $g_0$ and $g_1$ are concordant, and $\nabla_0$ and $\nabla_1$ can be joined by leafwise flat connections, then 
$i((g_0, \nabla_0),(g_1, \nabla_1)) = 0$.
\end{remark}

Next, we have a corollary of Theorem \ref{RelIndThm2}.  
\begin{coro}\label{cor6.4} Suppose $g_0, g_1, g_2 \in  \cM_{sc}^+$.
Then
$$
i(g_0,g_1) + i(g_1,g_2) \,\, =  \,\,  i(g_0,g_2), 
 \;\; \text{so} \;\; i(g_0,g_1) + i(g_1,g_2) +i(g_2,g_0) \,\, =  \,\, 0.
$$
\end{coro}

\begin{proof}  In the notation of Theorem \ref{RelIndThm2}, take $(M,F)$,  $(M',F')$ and $(\what{M}, \what{F})$ to be $(M_{\R}, F_{\R})$, $\cK = \cK' = M \times [0,1]$, $V_{\Phi} =  V_{\Phi}' =  M \times (-\infty, 0)$, and $V_+ = V_+' = M \times (1,\infty)$.
To compute $i(g_i,g_j)$ take 
$$
G_{i,j} = g_i + dt^2 \text{  for $t \in (-\infty,0]$, and  } G_{i,j} = g_j + dt^2  \text{  for  }t \in [1,\infty).
$$ 

For the first, we have
$$
i(g_0,g_1) - i(g_0,g_2)  =   {{\ch}}(\Ind_a(D_{\R}(G_{0,1}))) -  {{\ch}}( \Ind_a(D_{\R}(G_{0,2}))) =  
$$
$$
{{\ch}}(  \Ind_a(D_{\R}(G_{2,1}))= i(g_2,g_1)  =  - i(g_1,g_2).
$$
The second equality is from Theorem \ref{RelIndThm2}, where $\what{D}_L^{\what{E}} = D_{\R}(G_{2,1})$, $D_L^E = D_{\R}(G_{0,1})$, and $D_L^{E'} = D_{\R}(G_{0,2})$.

\ms

The second equation is then obvious, as $i(g_0,g_2) = -i(g_2,g_0)$.
\end{proof}

Now suppose that $M$ is the boundary of a compact manifold $W$ with a spin foliation $\what{F}$  which is transverse to $M$, and  which restricts to $F$ there.  Extend $\what{F}$ as above to  $W \cup_{M}  (M \times [0,\infty))$.   Given a metric $g$ of PSC on $F$, extend it to a complete leafwise metric $\what{g}$ on $\what{F}$  by making it $g + dt^2$ on $M \times [-\ep,\infty)$, where $M \times [-\ep,0]$ is a collar neighborhood  of $M \subset W$, and extending it arbitrarily over the rest of the interior of $W$.  
\begin{defi}\label{defig}
\hspace{0.5cm}$i(g,W) \,\, = \,\, {{\ch}}(\Ind_a(\what{D}_{\what{F}}))$.
\end{defi}
Note that Theorem \ref{classicalindex} and the proof of Theorem \ref{Phithm},  show that  $i(g,W)$ does not depend on the extension of $g$ over $W$.  {{It does however  depend on $W$ in general.}}

\ms 

In this situation, we have the following two corollaries of Theorem \ref{classicalindex}. 
\begin{coro}\label{cor6.6}  Suppose that $g_0, g_1 \in \cM_{sc}^+$.  Then 
$$
i(g_0,g_1) \; = \;  i(g_1,W) \; - \;  i(g_0,W), 
$$
as Haefliger classes.  In addition, if $\what{g}_0$ has PSC, then $i(g_0,W) = 0$. 
\end{coro}
The reader may wonder how the classes in the first equality can be compared, since they are on different manifolds.  This is explained below. 
\begin{proof}
Consider the following
\begin{itemize}
\item $(M_{\R},F_{\R})$ with the metric $G_{0,1}$ above, giving $i(g_0,g_1)$.
\item $M_0 = W_0 \cup_M  (M \times [0,\infty))$ with the metric $g_0 + dt^2$ on $M \times [0,\infty)$, and the metric $\what{g}_0$ on $W_0 = W$.  Take the opposite orientation on $M_0$ by reversing the orientations on $[0,\infty)$ and $W_0$, so this gives $-i(g_0,W)$. 
\item $M_1 = W_1 \cup_M  (M \times [0,\infty))$ with the metric $G_{0,1}$ restricted to $M \times [0,\infty)$, and the metric $\what{g}_0$ on $W_1 = W$.  As the metric on $M \times [1,\infty)$ is $g_1 \times dt^2$, this gives $i(g_1,W)$.  
\end{itemize} 
  
The meaning of the first equality is that representatives of the classes on $M_0 \, \dot\cup \, M_1 \ssm W_0  \, \dot\cup \,  W_1$ equal the representative on $M_{\R}$, while what remains on $W_0$ and $W_1$ cancel.  It is useful to have a picture of the situation.  The arrows indicate the orientations.

\begin{picture}(500,100)

\put(200,95){ $\line(0,-1){25}$}
\put(300,95){ $\line(0,-1){25}$}
\put(90,70){ $\line(1,0){320}$}
\put(90,95){ $\line(1,0){320}$}
\put(6,80){$i(g_0,g_1): \,\, M_{\R}$}
\put(95,80){$\cdots$}
\put(130,80){$g_0 + dt^2$}
\put(137,72){$\longrightarrow$}
\put(250,80){$g_t$}
\put(245,72){$\longrightarrow$}
\put(350,80){$g_1 + dt^2$}
\put(359,72){$\longrightarrow$}
\put(400,80){$\cdots$}
\put(186,60){$M \times\{0\}$}
\put(286,60){$M \times\{1\}$}

\put(150,55){ $\line(0,-1){25}$}
\put(200,55){ $\line(0,-1){25}$}
\put(300,55){ $\line(0,-1){25}$}
\put(150,30){ $\line(1,0){260}$}
\put(150,55){ $\line(1,0){260}$}
\put(6,40){$i(g_1,W): \,\, M_1$}
\put(160,40){$\what{g_0}$}
\put(170,31){$\longrightarrow$}
\put(180,40){$W_1$}
\put(250,40){$g_t$}
\put(245,31){$\longrightarrow$}
\put(350,40){$g_1 + dt^2$}
\put(359,31){$\longrightarrow$}
\put(400,40){$\cdots$}

\put(150,18){ $\line(0,-1){25}$}
\put(200,18){ $\line(0,-1){25}$}
\put(150,18){ $\line(1,0){260}$}
\put(150,-7){ $\line(1,0){260}$}
\put(-1,3){$-i(g_0,W): \,\, M_0$}
\put(160,4){$\what{g_0}$}
\put(180,4){$W_0$}
\put(170,-4){$\longleftarrow$}
\put(286,5){$g_0 + dt^2$}
\put(292,-4){$\longleftarrow$}
\put(400,5){$\cdots$}

\end{picture}

\vspace{0.5cm}

We may use the $\what{A}$-forms associated to the terms, since they are arbitrarily close to differential forms in the Haefliger classes.  We indicate them by $\what{A}(M_{\R})$, $\what{A}(M_0)$, and $\what{A}(M_1)$.  Then,
\begin{itemize} 

\item $\what{A}(M_1)$ restricted to $M_1 \ssm W_1$ equals $\what{A}(M_{\R})$ restricted to $M \times (0,\infty)$;
\item  $\what{A}(M_0)$ restricted to $M_0 \ssm W_0$ equals $\what{A}(M_{\R})$ restricted to $M \times (-\infty,0]$;
\item $\what{A}(M_0)$ restricted to $W_0$ cancels $\what{A}(M_1)$ restricted to $W_1$. 
\end{itemize} 
For the second statment, Propsition \ref{measPR2} gives that there is $\ep > 0$ so that $P_{[0,\ep]} = 0$.  Then Theorem \ref{HaefTheorem} gives  $i(g_0,W) = 0$.
\end{proof}

\begin{coro}\label{cor6.6bis} 
Suppose that $\what{g}_0$ has PSC, and that $g_1$ extends to $\what{g}_1$ with PSC over a compact manifold $\what{W}_1$ with the spin foliation $\what{F}_1$ extending $F$.  Set 
$$
X_{(0,1)} \;\ = \; W \cup_M  (M \times [0,1]) \cup_M \what{W}_1
$$
 with the  metric $\what{G}_{(0,1)}$ which is  $\what{g}_0$ on $W$,  $G_{0,1}$ on $M \times [0,1]$ and $\what{g}_1$ on $\what{W}_1$.
Denote the leafwise operator on the foliation $F_{(0,1)}$ of $X_{(0,1)} $ by $D_{(0,1)}$.  Then 
$$
i(g_0,g_1) \; = \;  \ch(\Ind_a(D_{(0,1)}))  \; = \; \int_{F_{(0,1)}} \hspace{-0.5cm} \what{A}(TF_{(0,1)}).
$$
\end{coro}

\begin{proof} 
For $i(g_0,g_1) =  \ch (\Ind_a(D_{(0,1)}))$, set $\what{M}_1 = \what{W}_1 \cup_M  (M \times [0,\infty))$ with the metric $g_1 + dt^2$ on $M \times [0,\infty)$, and the metric $\what{g}_1$ on $\what{W}_1$, so the metric has PSC everywhere and  $i(g_1,\what{W}_1) = 0$.  Then, we have,
$$
i(g_0,g_1)   \; = \; i(g_1,W) - i(g_0,W))  \; = \; i(g_1,W)   \; = \; i(g_1,W) - i(g_1,\what{W}_1) \; = \; 
{{\ch}} (\Ind_a(D_{(0,1)})).
$$
The first three equalities are obvious.  For the last, procede as in the first part, noting that 
\begin{itemize} 
\item $\what{A}(M_1)$ restricted to $W_1 \cup_M  (M \times [0,1])$ equals $\what{A}(X_{(0,1)})$ restricted to $W \cup_M  (M \times [0,1])$; 
\item  -$\what{A}(\what{M}_1)$ restricted to $\what{W}_1$ cancels $\what{A}(X_{(0,1)})$ restricted to $\what{W}_1$;
\item -$\what{A}(\what{M}_1)$restricted to $\what{M}_1 \ssm \what{W}_1$ cancels $\what{A}(M_1)$ restricted to $M_1 \ssm (M_1 \times (1,\infty))$.
\end{itemize}

Finally, the fact that $\dd {{\ch}} (\Ind_a(D_{(0,1)}))  = \int_{F_{(0,1)}} \hspace{-0.5cm} \what{A}(TF_{(0,1)})$ is a result from \cite{BH04}.
\end{proof}

\subsection{Some examples}

To finish, we construct a large collection of spin foliations whose space of leafwise PSC metrics has infinitely many path connected components.

\ms

Suppose we have the following data.
\begin{itemize}
\item A closed foliated manifold $(M,F)$, with $F$ spin and  $\dd \int_F \what{A}(TF)  \neq 0$ in $H_c^*(M/F)$. 
\item A closed manifold $S$ and a family $(g_i)$ of PSC metrics on it, and compact spin manifolds $X_i$ with boundary $S$ and metric $\what{g}_i$, which is $g_i \times dt^2$ in a neighborhood of  $S$, and $\what{g}_i$ also has PSC.   
Set 
$$
X_{(i,j)} \; = \; X_i  \cup (S \times [0,1]) \cup X_j,
$$
where the metric on $S \times [0,1]$ is $g_t \times dt^2$, and $g_t$ is a path of metrics from $g_i$ to $g_j$.
Assume further that $i(g_i ,g_j)$ is non-zero.  
\end{itemize}
\begin{proposition}\label{ExampleGL}
The foliated manifold $(M \times S,TF \times TS)$ has a family of PSC metrics $(\fg_i)$, so that for any $i \neq j$, $\fg_i$ and $\fg_j$ do not belong to  the same path component of the space of PSC metrics on $TF \times TS$. 
\end{proposition}
\begin{proof}  Since $M$ is compact, $F$  admits a  metric $g$ of bounded scalar curvature.
Set $\fg_i = g \times c_i g_i$, where $c_i \in (0,\infty)$ is such that $\fg_i$ has PSC.  
For the manifold $ M \times X_{(i,j)}$, with the foliation $F \times X_{(i,j)}$,  Corollary \ref{cor6.6bis} gives
 $$
i(\fg_i ,\fg_j) \, = \, \int_{F \times X_{(i,j)}} \hspace{-0.8cm} \what{A}(TF \times TX_{(i,j)}).
$$
 If $\fg_i$ and $\fg_j$ were in the same path component of the space of PSC metrics on $TF \times TS$, then we would have $i(\fg_i ,\fg_j) = 0$.  However, if $i \neq j$, then
$$
\int_{F \times X_{(i,j)}}  \hspace{-0.8cm} \what{A}(TF \times TX_{(i,j)}) =
\int_{F \times X_{(i,j)}} \hspace{-0.8cm}  \what{A}(TF)\what{A}(TX_{(i,j)}) =
 \int_F \what{A}(TF)\int_{X_{(i,j)}} \hspace{-0.5cm}\what{A}(TX_{(i,j)}) =   i(\fg_i ,\fg_j) \int_F \what{A}(TF) \neq 0.
$$
\end{proof}

Here are examples of this type.

\begin{example}\label{SLexample}

We adapt Example 1 of \cite{H78}.  In particular, let $G= SL_2\R \times \cdots \times SL_2\R$ (q copies) and $K= SO_2  \times \cdots  \times SO_2$ (q copies). $G$ acts naturally on $\R^{2q}\ssm \{0\}$ and is well known to contain subgroups $\Gam$ with $N = \Gam \backslash G/K$ compact, (in fact a product of $q$ surfaces of higher genus).   Set 
$$
M = \Gam \backslash G \times _K ((\R^{2q} \ssm \{0\})/\Z)  \simeq\Gam \backslash G \times _K (\S^{2q-1} \times \S^1),
$$
 where $n \in \Z$ acts on $\R^{2q} \ssm \{0\}$ by $n\cdot z = e^n z$. 
 
 $M$ has two transverse foliations, $F$ which is given by the fibers $\S^{2q-1} \times \S^1$ of the fibration $M \to N$, and a transverse foliation coming from the foliation $\tau$ of Example 1 of \cite{H78}.  More precisely, $\tau$ is defined on the vector bundle $\Gam \backslash G \times _K \R^{2q}$, and the zero section is a leaf of it.  In addition, the action of $\Z$ preserves $\tau$, fixing the zero section, so it descends to a foliation on $M$, also denoted $\tau$.
\end{example}
  
We work with $F$, noting that $TF$ is orientable and spin since $\R^{2q} - \{0\}$ has these structures and the actions of $K$ and $\Z$  preserve them.  It also happens to admits a metric with PSC, namely the product of the standard metrics on $\S^{2q-1}$ and $\S^1$, which is preserved by the action of $K$. {{The following proposition is proven in the appendix}}.  
\begin{prop}  \label{volFormOnN}
$\dd \int_F \what{A}(TF)$ is a nowhere zero $2q$ form on $N$.  In particular, there is a non-zero constant $C_q$ so that $\dd \int_N \int_F \what{A}(TF) = C_q \vol(N)$.
\end{prop}
Thus, $\dd \int_F \what{A}(TF) \neq 0$ in $H_c^*(M/F)$.  Note that this also shows that the Haefliger $\what{A}$ genus of $TF$, i.e. $\dd \left[ \int_F  \what{A}(TF) \right]^0 \in H_c^0(M/F)$, is zero, which is why we cannot use the results of \cite{BH21} here.
 
 \ms

In \cite{C88}, Carr constructs  examples of ``exotic" PSC metrics $\fg_i$, $i \in \Z_+$ on $\S^{4k-1}$, for $k > 1$, and compact Riemannian $4k$ dimensional spin manifolds $X_i$ with boundary $\S^{4k-1}$, so that the metric $\what{\fg}_i$ on $X_i$ is $\fg_i \times dt^2$ in a neighborhood of  $\S^{4k-1}$, and $\what{\fg}_i$ also has PSC.   
Set 
$$
X_{(i,j)} \; = \; X_i  \cup (\S^{4k-1} \times [0,1]) \cup X_j,
$$
where the metric on $\S^{4k-1} \times [0,1]$ is $\fg_t \times dt^2$, and $\fg_t$ is a path of metrics from $\fg_i$ to $\fg_j$.
These examples have the property that  $i(\fg_i ,\fg_j) \neq 0$. 
Thus we have all the elements required to apply Proposition \ref{ExampleGL}

\begin{remark}
Note that the calculations in the examples in \cite{H78} can be used  to provide examples associated to the groups  $G= SL_{2n_1}\R \times \cdots \times SL_{2n_r}\R$, and $K= SO_{2n_1}  \times \cdots  \times SO_{2n_r}$, and $G= SL_{2n_1}\R \times \cdots \times SL_{2n_r}\R \times \R$ and $K= SO_{2n_1}  \times \cdots  \times SO_{2n_r} \times \Z$.   We leave the details and further extensions to the reader.
\end{remark}

The next example is an easy corollary of  the Kreck-Stolz result from \cite{KS93}[Corollary 2.15].

\begin{proposition} Suppose that $(M,F)$ is a closed foliated manifold with $F$ spin.
Let $Y$ be a closed connected spin manifold of dimension $4k - 1 > 3$ with vanishing real Pontrjagin classes and such that $H^1(Y; \Z/2) = 0$.  If $Y$ admits a PSC metric, then  the foliated manifold $(M \times Y, TF\times TY)$ admits a sequence {{ $(\fg_i)$ of leafwise PSC metrics such that for any $i\neq j$, $\fg_i$ and $\fg_j$}}  are not in the same path component of PSC metrics on $TF \times TY$. 
\end{proposition}

Notice that if $Y$ is for instance simply connected, then it always admits a metric of PSC by \cite{St92}.

\begin{proof} In \cite{KS93}, Kreck and Stolz  produce  an infinite sequence $g_i$ of PSC metrics on $Y$ such that for any $i\neq j$, {the Gromov-Lawson invariant $i_{GL} (g_i, g_j)\neq 0$.  Note that $i_{GL}(g_i,g_j)$ is the difference of the dimensions of the positive and negative parts of the kernel of {$D_{\R}$ on the manifold $Y_{\R}$. 
 Thus,  there is a non-trivial $L^2$ element $\zeta$ in the kernel of $D_{\R}$.  On the foliated manifold  $(M \times Y, TF \times TY)$ there is the sequence of PSC metrics ($\fg_i$), where $\fg_i$ is as in Proposition \ref{ExampleGL}.}}  For $i \neq j$, these metrics are not in the same path component of leafwise PSC metrics.  For if they were, then the foliation $TF \times TY_{\R}$, would have PSC everywhere.  So, by Proposition \ref{measPR2},  there would not be any non-trivial $L^2$ elements in the kernel of $D_{\R}$. But this is patently false as $(0,\zeta)$ is such a non-trivial $L^2$ element. 
\end{proof}

\appendix

\section{Proof of Proposition \ref{volFormOnN}}  

We follow the proof of Theorem 5.4 in \cite{H78}.  Denote by $(x_1,y_1,.., x_q, y_q)$ the coordinates on $\R^{2q}$. Choose nonzero numbers $\lam_1,...\lam_q  \in \R$, and set
$$
X_{\lam} \, = \,  \sum_{i=1}^q \lam_i(x_i \pa/ \pa x_i + y_i \pa/\ \pa y_i).
$$
This vector field has an isolated singularity at the origin and it commutes with the actions of $K$ and  $\Z$ on $\R^{2q}  \ssm \{0\}$. Thus it induces a nowhere zero vector field also denoted $X_{\lam}$ on the bundle $M$.  

Denote by  $\ome_{\lam}$ the one-form on $\R^{2q} \ssm \{0\}$ defined by
$$
\ome_{\lam}  \, = \, \sum_{i=1}^q \frac{\lam_i}{x_i^2 + y_i^2}(x_i dx_i+ y_i dy_i).
$$
Note carefully that this is different from the $\ome_{\lam}$ of \cite{H78}.  This change is necessary so that $\ome_{\lam}$ is invariant under the action of $\Z$.  Note also that $d\ome_{\lam}  = 0$ still holds. 
The actions of $K$ and $\Z$ on $\R^{2q}  \ssm \{0\}$ preserve $\ome_{\lam}$, so it induces a one-form $\ome_{\lam}$ on M.

\ms

Let $S$ be the sphere bundle in $\what{M} = \Gam \backslash G \times _K ((\R^{2q} \ssm \{0\}))$ given by the image of 
$$
\{(g, (x_1,y_1,..., x_q,y_q)) \in G \times (\R^{2q} \ssm \{0\}) \, | \, \sum_{i=1}^q \lam_i(x_i^2 + y_i^2) = 1\}. 
$$
$S$ is invariant under $\Gam$ and $K$ so it is well defined.  Set 
$$
S_0 \, = \, S  \, = \, 0 \cdot S, \quad \text{and} \quad  S_1 = 1\cdot S.
$$
Note that the condition on $S_1$ is $ \sum_{i=1}^q \lam_i(x_i^2 + y_i^2) = e^2$, so its radius is $e$.  Then we may write,
$$
M \, = \,  \Gam \backslash G \times _K (\S^{2q-1} \times [1,e]),
$$
where we identify the boundary components, $S_0$ and $S_1$, on the right, and we may do our computations, as in \cite{H78}, using the coordinates on $G \times \S^{2q-1} \times (1,e)$.

\ms

Denote by  $\theta$ the unique basic connection (for the foliation $\tau$!) on $T(\S^{2q-1} \times (1,e))$, which is the normal bundle of $\tau$, whose covariant derivative $\nabla$  satisfies, for all $Y \in T(\S^{2q-1} \times (1,e))$, 
$$
\nabla_Y \pa/\pa x_i  \, =  \, \ome_{\lam} (Y)[X_{\lam}, \pa /\pa x_i], \;\; \text{and} \;\; \nabla_Y \pa/\pa y_i  \, =  \, \ome_{\lam} (Y)[X_{\lam}, \pa /\pa y_i]. 
$$
The proof in \cite{H78} works just as well here to show that $\theta$ is well defined.

\ms

The computation of the curvature $\Ome$ of $\nabla$ proceeds just as in \cite{H78}.  In particular,  we may assume that we have a neighborhood $U$ in $N$ whose inverse image in $M$ is of the form $U \times (\S^{2q-1} \times (1,e))$, and coordinates on it, so that the local form of  $\Ome$ with respect to the local basis $\pa/ \pa x_1,\pa/ \pa y_1, ... \pa/ \pa x_q, \pa/ \pa y_q$ of $T(\S^{2q-1} \times (1,e))$, is given by 
$$
\Ome_{2i-1}^{2i-1} \, = \, \Ome_{2i}^{2i}  \, = \, \lam_i d (\lam\delta),
$$
and all other terms are zero.
Here, for $i = 1,...,q$,
\begin{itemize}
\item  $\lam\delta \, = \, \lam_1\delta_1 + \cdots + \lam_q \delta_q$;
\item  $\delta_i = x_iy_i\ome_i + \frac{1}{2}(x_i^2 - y_i^2)\gam_i$;
\item $\ome_1,\gam_1,...,\ome_q,\gam_q$ is a basis of the one-forms on $U$   with $d\ome_i = -\ome_i \wedge \gam_i$ and $d \gam_i = 0$.
\end{itemize}
Recall that $\dd \what{A}(\xi_1,..., \xi_{2q}) = \prod_{j=1}^{2q}\frac{ \xi_j/2}{\sinh(\xi_j/2)}$.  We want to compute
$$
\int_F \what{A}(T(S^{2p-1} \times S^1)) \, = \,
\int_{S^{2q-1} \times S^1} \hspace{-1.0cm} \what{A}(T(S^{2p-1} \times S^1)) \, = \,
\int_{S^{2q-1} \times (1,e)}  \hspace{-1.0cm}\what{A}(\Ome).
$$
As $(d \delta_i)^3 = 0$ and $(d \delta_i)^2 = 2(x_i^2 + y_i)^2dx_i \wedge dy_i \wedge \ome_i \wedge \gam_i$, the only 
term of $\what{A}(\Omega)$ which will be non-zero when integrated over $F$ is, just as in \cite{H78},  
$$
\what{A}_{2q}(\Ome) \, = \, \what{A}_{2q}(\lam_1,\lam_1,...,\lam_q,\lam_q) (d (\lam \delta))^{2q} \, = \
\what{A}_{2q}(\lam_1,\lam_1,...,\lam_q,\lam_q)(2q)! \prod_{i=1}^q(\lam_i^2(x_i^2 + y_i^2) dx_i \wedge dy_i \wedge \ome_i \wedge \gam_i,        
$$
where $\what{A}_{2q}(\xi_1,...,\xi_{2q})$ is the term in $\what{A}(\xi_1,...,\xi_{2q})$ of degree 2$q$.  Thus,
$$
\int_{S^{2q-1} \times (1,e)} \hspace{-1.0cm} \what{A}_{2q}(\Ome) \, = \,
\what{A}_{2q}(\lam_1,\lam_1,...,\lam_q,\lam_q)\left[(2q)! \int_{S^{2q-1} \times (1,e)} 
\prod_{i=1}^q(\lam_i^2(x_i^2 + y_i^2) dx_i \wedge dy_i)\right]    \prod_{i=1}^q(\ome_i \wedge \gam_i)\, = \,
$$
$$
\frac{\pi^q(e^{4q}-1)\what{A}_{2q}(\lam_1,\lam_1,...,\lam_q,\lam_q)}{(\lam_1 \cdots \lam_q)^2}    \prod_{i=1}^q(\ome_i \wedge \gam_i),
$$ 
by Lemma 5.8 of \cite{H78}, which is a nowhere zero $2q$ form on $N$.  Note that $\what{A}_{2q}(\lam_1,\lam_1,...,\lam_q,\lam_q)$ is a non-zero constant times $(\lam_1 \cdots \lam_q)^2$.  Thus, there is a non-zero constant $C_q$ so that 
 
\hspace{4.0cm} 
$
\dd \int_N \int_{S^{2q-1} \times (1,e)} \hspace{-1.0cm} \what{A}(\Ome) =  \dd \int_N \int_F \what{A}(TF) = C_q \vol(N).
$


\begin{thebibliography}{MMM11t}

\bibitem [A76]{A76} M. F. Atiyah,  {\em Elliptic operators, discrete groups and Von Neumann algebras}, Asterisque, Societe Math. de France  {\bf 32-33} (1976) 43--72.

\bibitem[AS68III]{ASIII} M. F. Atiyah and I. M. Singer,    {\em  The index of elliptic operators III},  Ann.  of Math. {\bf 87} (1968) 546--604.

\bibitem[Be20]{Be20} M-T. Benameur,  {\em  The relative L2 index theorem for Galois coverings.}  Ann. K-Theory 6 (2021), no. 3, 503-541. 

\bibitem[BH04]{BH04} M-T. Benameur and J. L. Heitsch,  {\em  Index theory and Non-Commutative Geometry I. Higher Families Index Theory}, K-Theory  {\bf 33} (2004) 151--183.   
{\em Corrigendum},   ibid  {\bf 36} (2005) 397--402.

\bibitem[BH08]{BH08} M-T. Benameur and J. L. Heitsch,  {\em Index theory and noncommutative geometry II.  Dirac operators and index bundles},  Journal of K-Theory {\bf 1} (2008) 305--356.

\bibitem[BHW14]{BHW14} M-T. Benameur, J. L. Heitsch \& Charlotte Wahl,  {\em  An interesting example for spectral invariants}, Journal of K-Theory {\bf 13} (2014) 305--311. 

\bibitem[BH18]{BH18} M-T. Benameur and J. L. Heitsch,  {\em Transverse non-commutative geometry of foliations},  
J. Geom. Phys. {\bf 134} (2018) 161--194. 

\bibitem[BH21]{BH21} M-T. Benameur and J. L. Heitsch,  {\em Dirac operators on foliations with invariant transverse measures},  J. Funct. Anal. 284 (2023), no. 2, Paper No. 109742, 52 pp.

\bibitem[BH23]{BH23} M-T. Benameur and J. L. Heitsch,  {\em Higher relative index theorems for foliations, 
a superconnection approach},  in preparation.

 \bibitem[B86]{B86} J.-M. Bismut, {\em The Atiyah-Singer index theorem for families of Dirac operators: two heat equation proofs.}  Invent. Math.  {\bf 83} (1986), 91--151.

\bibitem[C88] {C88} R. Carr, {\em Construction of manifolds of positive scalar curvature}, Trans. Amer. MAth. Soc. {\bf 307}, (1988) 63--74.

\bibitem[Ch73]{Ch73} P. R. Chernoff,  {\em Essential self-adjointness of powers of generators of hyperbolic equations},  
J. Functional Analysis {\bf 12} (1973), 401--414. 

\bibitem[C79]{C79} A. Connes,  {\em Sur la th\'eorie non commutative de l'int\'egration.} (French) Alg\`ebres d'op\'erateurs (S\'em., Les Plans-sur-Bex, 1978), pp. 19-143, Lecture Notes in Math., 725, Springer, Berlin, 1979.

\bibitem[C81]{C81} A.~Connes, {\em A survey of foliations
and operator algebras}.  Operator algebras and applications, Part I,
Proc. Sympos.  Pure Math {\bf 38}, Amer.  Math.  Soc. (1982) 521--628.

\bibitem[CS84]{CS84} A. Connes and G. Skandalis, {\em The longitudinal index theorem for foliations.} Publ. Res. Inst. Math. Sci. 20 (1984), no. 6, 1139-1183. 

\bibitem[C94]{C94} A. Connes, 
Noncommutative Geometry, 
Academic Press, New York - London, 1994.

\bibitem[Do87]{Do87} H. Donnelly, {\em Essential spectrum and heat kernel.} J. Funct. Anal. 75 (1987), no. 2, 362-381. 

\bibitem[GL83]{GL83} M. Gromov and H. B. Lawson, Jr.,  {\em Positive scalar curvature and The Dirac operator on complete Riemannian manifolds},  Publ. Math. I.H.E.S. {\bf 58} (1983) 295--408.

\bibitem[Ha80]{Ha80} A. Haefliger, {\em Some remarks on foliations with minimal leaves}, J. Diff. Geo. {\bf 15}  (1980) 269--284.

\bibitem [H78]{H78} J. L. Heitsch,  {\em  Independent variation of secondary classes},  Ann. of Math. {\bf 108} (1978) 421--460. 

\bibitem [H95]{H95} J. L. Heitsch,  {\em  Bismut superconnections and the {C}hern character for {D}irac
 operators on foliated manifolds}, K-Theory {\bf 9} (1995) 507--528.

\bibitem[HL90]{HL90} J. L. Heitsch and C. Lazarov, {\em A Lefschetz theorem for foliated manifolds},  Topology {\bf 29} (1990) 127--162.

\bibitem[HL99]{HL99} J. L. Heitsch and C. Lazarov, {\em A general families index theorem}, K-Theory {\bf 18} (1999) 181--202.

\bibitem[KS93] {KS93}M. Kreck and S. Stolz, {\em Nonconnected moduli spaces of positive sectional curvature metrics}, J. of AMS {\bf 6} (1993)  825--850.

\bibitem[LM89]{LM89} H.  B. Lawson, Jr. and  M.-L. Michelson, {\em Spin geometry}, Princeton Math.  Series {\bf 38}, Princeton, 1989.

\bibitem[NWX96]{NistorWeinsteinXu} V. Nistor, A. Weinstein and P. Xu,
{\em Pseudodifferential operators on differential groupoids},  Pacific
J. Math.  {\bf 189} (1999) 117--152.

\bibitem[R87]{Roe87} J. Roe,   {\em Finite propagation speed and Connes' foliation algebra}, Math. Proc. Cambridge Philos. Soc. {\bf 102} (1987) 459--466. 

\bibitem[Sh92]{Shubin} M.A. Shubin,   {\em Spectral theory of elliptic operators on non-compact manifolds}, M\'{e}thodes Semi-classiques, Vol. 1 (Nantes, 1991),  Ast\'erisque  {\bf 207} (1992) 35--108.

\bibitem[St92]{St92} S. Stolz, {\em Simply connected manifolds of positive scalar curvature}, Ann. of Math. (2) {\bf 136}
(1992) 511--540.

\bibitem[Vi67]{Vi67} E. Vesentini, {\em Lectures on Levi convexity of complex manifolds and cohomology vanishing theorems.} Notes by M. S. Raghunathan
Tata Inst. Fundam. Res. Lect. Math., No. 39
Tata Institute of Fundamental Research, Bombay, 1967. iii+131+vi pp.
(1992) 511-540.

\end{thebibliography}
\end{document}